\input amstex
\documentstyle{amsppt}

\magnification=\magstep1
\NoBlackBoxes
\NoRunningHeads
\topmatter
\title Actions of symbolic dynamical systems on $C^*$-algebras II.\\
Simplicity of $C^*$-symbolic crossed products and some examples 
\endtitle
\author Kengo Matsumoto
\endauthor
\affil Department of Mathematical Sciences \\
 Yokohama City University\\
  Seto 22-2, Kanazawa-ku, Yokohama 236-0027, JAPAN
\endaffil
\abstract
We have introduced a notion of $C^*$-symbolic dynamical system
in [K. Matsumoto: Actions of symbolic dynamical systems on $C^*$-algebras,
 to appear in J. Reine Angew. Math.],
 that is a finite family of endomorphisms of a $C^*$-algebra  
with some conditions.
The endomorphisms are indexed by symbols and  yield both a subshift and a $C^*$-algebra of a Hilbert $C^*$-bimodule. 
The associated $C^*$-algebra with the $C^*$-symbolic dynamical system 
 is regarded as a crossed product by the subshift.
We will study a simplicity condition of the 
$C^*$-algebras of the $C^*$-symbolic dynamical systems.
Some examples 
such as irrational rotation Cuntz-Krieger algebras will be studied.
 \endabstract

\endthanks

\endtopmatter


\def\Ker{{{\operatorname{Ker}}}}
\def\id{{{\operatorname{id}}}}
\def\Ad{{{\operatorname{Ad}}}}
\def\Sp{{{\operatorname{Sp}}}}
\def\End{{{\operatorname{End}}}}
\def\Aut{{{\operatorname{Aut}}}}
\def\min{{{\operatorname{min}}}}

\def\dist{{{\operatorname{dist}}}}

\def\A{{ {\Cal A} }}
\def\B{{ {\Cal B} }}

\def\I{{ {\Cal I} }}
\def\E{{ {\Cal E} }}
\def\Zp{{ {\Bbb Z}_+ }}
\def\ARL{{{\Cal A}\rtimes_{\rho}\Lambda}}

\def\FN{{{\Cal F}_{\theta_1,\dots,\theta_N}}}
\def\FNk{{{\Cal F}^k_{\theta_1,\dots,\theta_N}}}
\def\FR{{ {\Cal F}_\rho}}
\def\FRK{{ {\Cal F}_\rho^k}}
\def\FRS{{{\Cal F}_{\rho^{\Sigma\otimes}}}}
\def\FRSK{{{\Cal F}_{\rho^{\Sigma\otimes}}^k}}
\def\BMA{{{\B\otimes\A}}}
\def\BMADS{{{(\BMA, \rho^{\Sigma\otimes},\Sigma)}}}
\def\CMA{{{{\Bbb C}\otimes\A}}}
\def\BMARL{{{(\BMA)\rtimes_{\rho^{\Sigma\otimes}}\Lambda}}}
\def\CMARL{{{({\Bbb C}\otimes\A)\rtimes_{\rho^{\Sigma\otimes}}\Lambda}}}
\def\DBMA{{{{\Cal D}_{\rho^{\Sigma\otimes}}}}}
\def\DA{{{{\Cal D}_{\rho}}}}
\def\DBMAk{{{{\Cal D}_{\rho^{\Sigma\otimes}}^k}}}
\def\DAk{{{{\Cal D}_{\rho}^k}}}
\def\OGN{{{\Cal O}_{{\Cal G}, \theta_1,\dots,\theta_N}}}
\def\ON{{{\Cal O}_{\theta_1,\dots,\theta_N}}}
\def\FGN{{{\Cal F}_{{\Cal G},\theta,\dots,\theta_N}}}


\heading 1. Introduction
\endheading
In [CK], 
 J. Cuntz and W. Krieger have founded a close relationship
 between symbolic dynamics and $C^*$-algebras (cf.[C], [C2]).
They constructed purely infinite simple $C^*$-algebras 
from irreducible topological Markov shifts.
The $C^*$-algebras are called  Cuntz-Krieger algebras.

In [Ma], the author introduced a notion of $\lambda$-graph system, whose matrix version is called symbolic matrix system. 
 A $\lambda$-graph system is a generalization of finite labeled graph and
 presents a subshift. 
He constructed $C^*$-algebras from $\lambda$-graph systems [Ma2]
as a generalization of the above Cuntz-Krieger algebras. 
A $\lambda$-graph system gives rise to a finite family 
$\{ \rho_\alpha \}_{\alpha \in \Sigma}$ of endomorphisms of 
a  unital commutative AF-$C^*$-algbera $\A_{\frak L}$ 
with some conditions stated below. 
A $C^*$-{\it symbolic dynamical system},\ 
introduced in [Ma6], is a generalization of $\lambda$-graph system.
It is a finite family $\{ \rho_{\alpha} \}_{\alpha \in \Sigma}$ 
of endomorphisms of a unital $C^*$-algebra 
$\Cal A$ 
such that the closed ideal generated by $\rho_\alpha(1), \alpha \in \Sigma$
coincides with $\A$.
A  finite labeled graph gives rise to a $C^*$-symbolic dynamical system 
$({\Cal A},\rho, \Sigma)$ such that $\Cal A = {\Bbb C}^N$ for some 
$N \in \Bbb N$.
Conversely, if $\Cal A = {\Bbb C}^N $, 
the $C^*$-symbolic dynamical system comes from a  finite labeled graph.
A  $\lambda$-graph system $\frak L$ 
gives rise to a $C^*$-symbolic dynamical system $({\Cal A},\rho, \Sigma)$ 
such that $\Cal A$ is $C(\Omega_{\frak L})$ for some compact Hausdorff space 
$\Omega_{\frak L}$ with 
$\dim \Omega_{\frak L} = 0$.
Conversely, if $\Cal A$ is $C(X)$ for a compact Hausdorff space 
$X$ with 
$\dim X = 0$, 
the $C^*$-symbolic dynamical system comes from a  $\lambda$-graph system.

A $C^*$-symbolic dynamical system $({\Cal A},\rho, \Sigma)$
yields a nontrivial subshift 
$\Lambda_{({\Cal A}, \rho, \Sigma)}$, 
that we will denote by $\Lambda_\rho$,  
over $\Sigma$
and a Hilbert $C^*$-right $\Cal A$-module 
$(\phi_\rho, {\Cal H}_{\Cal A}^{\rho}, \{u_\alpha\}_{\alpha \in \Sigma})$
that has an orthogonal finite basis 
$\{u_{\alpha}\}_{\alpha\in \Sigma}$
and a unital faithful diagonal left action 
$\phi_{\rho} :{\Cal A} \rightarrow L({\Cal H}_{\Cal A}^{\rho})$.
It is called a Hilbert $C^*$-{\it symbolic bimodule over}\ $\Cal A$,
and written as
$(\phi_{\rho},{\Cal H}_{\Cal A}^{\rho}, 
\{u_{\alpha}\}_{\alpha\in \Sigma})$.                                           
By using general construction of $C^*$-algebras 
 from Hilbert $C^*$-bimodules established by M. Pimsner [Pim]
 (cf. [Ka]),
 the author has introduced  
a $C^*$-algebra
denoted  by 
${\Cal A}\rtimes_{\rho}\Lambda$
from the Hilbert $C^*$-symbolic bimodule 
$(\phi_{\rho},{\Cal H}_{\Cal A}^{\rho}, 
\{u_{\alpha}\}_{\alpha\in \Sigma})$, 
where 
$\Lambda$ is the subshift $\Lambda_{\rho}$
associated with 
$({\Cal A}, \rho, \Sigma)$.
We call the algebra 
${\Cal A}\rtimes_{\rho}\Lambda$ 
{\it the}\ $C^*$-{\it symbolic crossed product of}\ $\Cal A$ {\it by the subshift}\ $\Lambda$.
If $\Cal A =\Bbb C$, 
the subshift $\Lambda$ is the full shift $\Sigma^{\Bbb Z}$,
and the $C^*$-algebra $\ARL$ 
is the Cuntz algebra ${\Cal O}_{|\Sigma|}$ of order  $|\Sigma|$. 
If $\Cal A = C(X)$ with $\dim X=0$, 
there uniquely exists a $\lambda$-graph system $\frak L$ up to equivalence
 such that the subshift $\Lambda$ 
is presented by $\frak L$ 
and the $C^*$-algebra $\ARL$ is the $C^*$-algebra ${\Cal O}_{\frak L}$
associated with the $\lambda$-graph system $\frak L$.
Conversely, for any subshift, that is presented by a $\lambda$-graph system $\frak L$,  
there exists a $C^*$-symbolic dynamical system  
$({\Cal A}, \rho, \Sigma)$ such that $\Lambda_{\rho}$
is the subshift 
presented by $\frak L$, 
the algebra $\Cal A$ is $C(\Omega_{\frak L})$  with 
$\dim \Omega_{\frak L} = 0$, 
and the algebra $\ARL$ is the $C^*$-algebra ${\Cal O}_{\frak L}$
associated with $\frak L$ ([Ma6]).
If in particular, $\Cal A = {\Bbb C}^n$, 
the subshift $\Lambda$ is a sofic shift
and 
$\ARL$ is a Cuntz-Krieger algebra.

In this paper, a condition called (I) on $(\A,\rho,\Sigma)$ 
is introduced
as a generalization of condition (I) on the finite matrices of 
Cuntz-Krieger [CK] and on the $\lambda$-graph systems [Ma2].
Under the assumption that $(\A,\rho,\Sigma)$ satisfies condition (I),
the simplicity conditions of 
the algebra $\ARL$ is discussed in Section 3.
We further study ideal structure of $\ARL$ 
from the view point of quotients of the $C^*$-symbolic dynamical systems
in Section 4.
Related discussions have been studied in Kajiwara-Pinzari-Watatani's paper [KPW] for the $C^*$-algebras of Hilbert $C^*$-bimodules 
(cf. [Kat], [MS], [Tom], etc.). 
They have studied simplicity condition and ideal structure of the $C^*$-algebras of Hilbert $C^*$-bimodules
in terms of the language of the Hilbert $C^*$-bimodules.
Our approach 
to study the algebras $\ARL$ is 
from the view point of $C^*$-symbolic dynamical systems,
that is differnt from theirs.
In Section 5, we will study pure infiniteness of the algebras $\ARL$. 
To obtain rich examples of the algebras $\ARL$, 
we will in Section 6 construct 
$C^*$-symbolic dynamical systems from a finite family of automorphisms
$\alpha_i \in \Aut(\B), i=1,\dots N$ on a unital $C^*$-algebra $\B$ 
and a $C^*$-symbolic dynamical systems $(\A,\rho,\Sigma)$
with
$
\Sigma = \{\alpha_1,\dots, \alpha_N\}.
$
The $C^*$-symbolic dynamical system is denoted by
$\BMADS$ that is the tensor product between two $C^*$-symbolic dynamical systems$(\B,\alpha,\Sigma)$ and $(\A,\rho,\Sigma)$.
 As examples of $C^*$-symbolic crossed products, 
 continuous analogue of Cuntz-Krieger algebras
 called irrational rotation Cuntz-Krieger algebras denoted by $\OGN$
 and    
irrational rotation Cuntz algebras denoted by $\ON$ are studied in Sections 8 and 9.
They belongs to the class of the $C^*$-algebras of continuous graphs 
by V. Deaconu ([De],[De2]).
The fixed point algebras $\FGN$ of $\OGN$ under gauge actions are no longer AF-algebras. 
They are A${\Bbb T}$-algebras.
In particular, the fixed point algebras $\FN$
of $\ON$ under gauge actions
are simple A${\Bbb T}$-algebras of real rank zero with unique tracial state
if and only if differnce of rotation angles $\theta_i-\theta_j$  is irrational for some $i,j=1,\dots,N$
(Theorem 9.4).  

 Throughout this paper, we denote by 
 $\Zp$ and by ${\Bbb N}$
 the set of nonnegative integers and 
 the set of positive integers respectively. 
A homomorphism and an isomorphism between  $C^*$-algebras
mean a $*$-homomorphism and a $*$-isomorphism respectively.
An ideal of a $C^*$-algebra means a closed two sided $*$-ideal.

\heading 2. $C^*$-symbolic dynamical systems and their crossed products
\endheading

Let ${\Cal A}$ be a unital $C^*$-algebra.
In what follows,
an endomorphism of $\Cal A$ means 
a $*$-endomorphism of $\Cal A$ that does not necessarily preserve the unit
$1_\A$ of 
$\Cal A$.
The unit $1_\A$ is denoted by $1$ unless we specify.
We denote by $\End({\Cal A})$ the set of all endomorphisms of ${\Cal A}$. 
Let $\Sigma$ be a finite set.
A finite family of endomorphisms 
$
\rho_\alpha \in \End({\Cal A}), \alpha \in \Sigma
$    
is said to be {\it essential}\ if
$\rho_{\alpha}(1) \ne 0$ for all $\alpha \in \Sigma$
 and   
the closed ideal generated by $\rho_\alpha(1), \alpha \in \Sigma$
coincides with $\A$.
It is said to be {\it faithful}\ if for any nonzero $x \in \Cal A$ there exists a symbol $\alpha\in \Sigma$ such that $\rho_{\alpha}(x) \ne 0$.
We note that 
 $\{ \rho_\alpha \}_{\alpha \in \Sigma}$ is faithful if and only if
the homomorphism
 $
 \xi_\rho : a \in \A \longrightarrow
  [ \rho_\alpha (a) ]_{\alpha \in \Sigma} \in \oplus_{\alpha \in \Sigma} \A
 $ 
is injective.

\noindent
{\bf Definition ([Ma6]).} A $C^*$-{\it symbolic dynamical system}\ 
is a triplet $({\Cal A}, \rho, \Sigma)$ 
consisting of a unital $C^*$-algebra $\A$ 
and 
an essential and faithful finite family of endomorphisms $\rho_{\alpha}$ of ${\Cal A}$ indexed by $\alpha\in \Sigma$.

Two $C^*$-symbolic dynamical systems $({\Cal A}, \rho, \Sigma)$
and
$({\Cal A}', \rho', \Sigma')$ 
are said to be isomorphic if there exist   
an isomorphism $\Phi :\Cal A \rightarrow {\Cal A}'$
and a bijection $\pi :\Sigma \rightarrow \Sigma'$ 
such that
$\Phi \circ \rho_{\alpha} = \rho'_{\pi(\alpha)}\circ \Phi$
for all 
$\alpha \in \Sigma$.
A $C^*$-symbolic dynamical system $({\Cal A}, \rho, \Sigma)$
 yields  a subshift $\Lambda_{({\Cal A}, \rho, \Sigma)}$
over $\Sigma$ such that a word $\alpha_1\cdots\alpha_k$ of $\Sigma$ is admissible for $\Lambda_{({\Cal A}, \rho, \Sigma)}$
 if and only if 
$(\rho_{\alpha_k}\circ \cdots\circ \rho_{\alpha_1})(1) \ne 0$
([Ma6;Proposition 2.1]).
The subshift $\Lambda_{({\Cal A}, \rho, \Sigma)}$
will be denoted by 
$\Lambda_{\rho}$ or simply by 
$\Lambda$
in this paper.

Let $\Cal G=(G,\lambda)$ be a left-resolving finite labeled graph 
with underlying finite directed graph $G =(V,E)$ and labeling map 
$\lambda:E \rightarrow \Sigma$ (see [LM; p.76]).
Denote by $v_1,\dots, v_N$ the vertex set $V$.
Assume that every vertex has both an incoming edge and an outgoing edge. 
Consider the $N$-dimensional commutative $C^*$-algebra 
$\Cal A_{\Cal G} = {\Bbb C}E_1\oplus\cdots \oplus {\Bbb C}E_N$ 
where each minimal projection $E_i$ 
corresponds to the vertex $v_i$ 
for
$ i=1,\dots,N$.
Define an $N\times N$-matrix for $\alpha \in \Sigma$
by 
$$
A^{\Cal G}(i,\alpha,j) =
\cases
1 & \text{ if there exists an edge } e \text{ from } v_i
\text{ to } v_j \text{ with } \lambda(e) = \alpha,\\
0 & \text{ otherwise}
\endcases      \tag 2.1
$$ for $i,j = 1,\dots,N$.
We set
$
\rho^{\Cal G}_{\alpha}(E_i) = \sum_{j=1}^N A^{\Cal G}(i,\alpha,j)E_j 
$
for
$ i= 1,\dots,N, \alpha \in \Sigma.
$ 
Then  
$\rho^{\Cal G}_{\alpha}, \alpha \in \Sigma$
define endomorphisms of ${\Cal A}_{\Cal G}$ 
such that 
$({\Cal A}_{\Cal G}, \rho^{\Cal G}, \Sigma)$ 
is a $C^*$-symbolic dynamical system 
such that the algebra ${\Cal A}_{\Cal G}$ is ${\Bbb C}^N$,
and the subshift 
$\Lambda_{\rho^{\Cal G}}$
is the sofic shift
$\Lambda_{\Cal G}$ 
presented by $\Cal G$.
Conversely,
for a $C^*$-symbolic dynamical system  
$({\Cal A}, \rho, \Sigma)$,
if $\Cal A $ is ${\Bbb C}^N$, 
there exists a left-resolving labeled graph 
$\Cal G$ such that $\Cal A = {\Cal A}_{\Cal G}$ 
and  
$\Lambda_{\rho} =\Lambda_{\Cal G} $
the sofic shift
presented by $\Cal G$ ([Ma6;Proposition 2.2]).

More generally let $\frak L$ be a $\lambda$-graph system 
$ (V,E,\lambda,\iota)$ over $\Sigma$ (see [Ma]).
Its vertex set $V$ is $\cup_{l=0}^{\infty}V_l$. 
We equip $V_l$ with discrete topology.
We denote by 
$\Omega_{\frak L}$ 
the 
compact Hausdorff space 
with $\dim \Omega_{\frak L} =0$
of the projective limit
$
V_0 \overset{\iota}\to{\leftarrow} V_1 \overset{\iota}\to{\leftarrow} 
V_2 \overset{\iota}\to{\leftarrow} \cdots,
$
as in [Ma2;Section 2].
The algebra $C(V_l)$ of all continuous functions on $V_l$,
denoted by ${\Cal A}_{{\frak L},l}$,
is the direct sum
$
{\Cal A}_{{\frak L},l} = {\Bbb C}E_1^l \oplus  
\cdots \oplus {\Bbb C}E_{m(l)}^l
$
where 
each
minimal projection
$E_i^l$
corresponds to
the vertex 
$v_i^l$ 
for 
$i=1,\dots,m(l)$.   
Let ${\Cal A}_{\frak L}$ be the commutative $C^*$-algebra
$
C(\Omega_{\frak L})
 =
\lim_{l \to \infty} \{ \iota_*:
{\Cal A}_{{\frak L},l}\rightarrow
{\Cal A}_{{\frak L},l+1} \}.
$
Let $A_{l,l+1}, l\in \Zp$ be the matrices defined in [Ma2; Theorem A].
For a symbol $\alpha \in \Sigma$ we set
$$
\rho_{\alpha}^{\frak L}(E_i^l) = \sum_{j=1}^{m(l+1)}A_{l,l+1}
(i,\alpha,j)E_j^{l+1}\qquad
\text{ for } i=1,2,\dots,m(l), \tag 2.2 
$$
so that 
$\rho_{\alpha}^{\frak L}$ defines an endomorphism of ${\Cal A}_{\frak L}$.
We have  a 
$C^*$-symbolic dynamical system 
$({\Cal A}_{\frak L}, \rho^{\frak L}, \Sigma)$
such that
 the $C^*$-algebra ${\Cal A}_{\frak L}$ is $C(\Omega_{\frak L})$ 
 with $\dim \Omega_{\frak L} =0$, 
and the subshift 
$\Lambda_{\rho^{\frak L}}$
coincides with the subshift $\Lambda_{\frak L}$ presented by $\frak L$. 
Conversely, for a 
$C^*$-symbolic dynamical system
 $({\Cal A}, \rho, \Sigma)$,
 if the algebra ${\Cal A}$ is $C(X)$ with $\dim X =0$,
there exists a $\lambda$-graph system $\frak L$ over $\Sigma$
such that the associated $C^*$-symbolic dynamical system 
$({\Cal A}_{\frak L}, \rho^{\frak L}, \Sigma)$
is isomorphic to 
$({\Cal A}, \rho, \Sigma)$ ([Ma6;Theorem 2.4]).

Let ${\frak L}$ and ${\frak L}'$ 
be predecessor-separated $\lambda$-graph systems over 
$\Sigma$ and $\Sigma'$ respectively. 
Then 
$({\Cal A}_{\frak L}, \rho^{\frak L}, \Sigma)$ is isomorphic to 
$({\Cal A}_{{\frak L}'}, \rho^{{\frak L}'}, \Sigma')$ if and only if
${\frak L}$ and ${\frak L}'$  are equivalent.  
In this case, the presented subshifts $\Lambda_{\frak L}$ 
and $\Lambda_{{\frak L}'}$ are identified through a symbolic conjugacy.
Hence the  equivalence classes of the $\lambda$-graph systems are identified with the isomorphism classes of the $C^*$-symbolic dynamical systems of the commutative AF-algebras.

We say that {\it a subshift}\ $\Lambda$ {\it acts on a}\ $C^*$-{\it algebra}\ 
${\Cal A}$ if there exists a $C^*$-symbolic dynamical system 
$({\Cal A}, \rho, \Sigma)$ 
such that the associated subshift $\Lambda_{\rho}$
is $\Lambda$.
For a $C^*$-symbolic dynamical system
$({\Cal A},\rho,\Sigma)$, 
we have a  Hilbert $C^*$-bimodule 
$(\phi_{\rho},{\Cal H}_{\Cal A}^{\rho}, \{u_{\alpha}\}_{\alpha\in \Sigma})$
called a Hilbert $C^*$-symbolic bimodule ([Ma6]).
We then have a $C^*$-algebra 
by using the Pimsner's general construction of 
$C^*$-algebras from  Hilbert $C^*$-bimodules
[Pim] 
(cf. [Ka], see also [KPW], [KW], [Kat], [MS], [PWY], [Sch] etc.).
We denote the $C^*$-algebra by ${\Cal A}\rtimes_{\rho}\Lambda$, where 
$\Lambda$ is the subshift $\Lambda_{\rho}$
associated with
$({\Cal A}, \rho, \Sigma)$.  
We call the algebra 
${\Cal A}\rtimes_{\rho}\Lambda$ 
{\it the}\ $C^*$-{\it symbolic crossed product of}\ ${\Cal A}$ {\it by the subshift}\ 
$\Lambda$. 
\proclaim{Proposition 2.1([Ma6;Proposition 4.1])}
The $C^*$-symbolic crossed product ${\Cal A}\rtimes_{\rho}\Lambda$
is the universal $C^*$-algebra
$C^*(x, S_{\alpha}; x \in \A, \alpha \in \Sigma)$  
generated by 
$x \in \A$ 
and partial isometries $S_{\alpha}, \alpha \in \Sigma$
subject to the following relations called $(\rho)$:
$$
\sum_{\beta \in \Sigma}S_{\beta}S_{\beta}^* =1,\qquad
S_\alpha^* x S_\alpha = \rho_\alpha(x),
\qquad
x S_\alpha S_\alpha^* =  S_\alpha S_\alpha^*  x  
$$
for all $x \in \Cal A$ and $\alpha \in \Sigma.$
Furthermore for  $\alpha_1,\dots,\alpha_k \in \Sigma$,  
a word $(\alpha_1,\dots,\alpha_k)$ is admissible for the subshift
$\Lambda$ 
if and only if
$S_{\alpha_1}\cdots S_{\alpha_k} \ne 0.$
\endproclaim
Assume that ${\Cal A}$ is commutative. 
Then we know ([Ma6;Theorem 4.2])
\roster
\item"(i)" 
If $\Cal A =\Bbb C$, 
the subshift $\Lambda$ is the full shift $\Sigma^{\Bbb Z}$,
and the $C^*$-algebra $\ARL$ is the Cuntz algebra ${\Cal O}_{|\Sigma|}$ of order  $|\Sigma|$. 
\item"(ii)"
If $\Cal A = {\Bbb C}^N$ for some $N \in \Bbb N$, 
the subshift $\Lambda$ 
is a sofic shift $\Lambda_{\Cal G}$ 
presented by a left-resolving labeled graph $\Cal G$,
and the $C^*$-algebra $\ARL$ is a Cuntz-Krieger algebra ${\Cal O}_{\Cal G}$ associated with the labeled graph.
Conversely, 
for any sofic shift $\Lambda_{\Cal G}$, 
that is presented by a left-resolving labeled graph 
$\Cal G$,  there exists a $C^*$-symbolic dynamical system  
$({\Cal A}, \rho, \Sigma)$ 
such that the associated subshift is the sofic shift $\Lambda_{\Cal G}$, 
the algebra $\Cal A$
is ${\Bbb C}^N$ for some $N \in \Bbb N$,   
and the $C^*$-algebra $\ARL$ is the Cuntz-Krieger algebra ${\Cal O}_{\Cal G}$ associated with the labeled graph ${\Cal G}$.
\item"(iii)"
If $\Cal A = C(X)$ with $\dim X =0$, 
there uniquely exists a $\lambda$-graph system $\frak L$ up to equivalence
 such that the subshift $\Lambda$ 
is presented by $\frak L$ 
and the $C^*$-algebra $\ARL$ is the $C^*$-algebra ${\Cal O}_{\frak L}$
associated with the $\lambda$-graph system $\frak L$.
Conversely, 
for any subshift
$\Lambda_{\frak L}$,
that is presented by a left-resolving $\lambda$-graph system $\frak L$,  
there exists a $C^*$-symbolic dynamical system  
$({\Cal A}, \rho, \Sigma)$ such that the associated subshift is the subshift 
$\Lambda_{\frak L}$, 
the algebra ${\Cal A}$ is $C(\Omega_{\frak L})$ with 
$\dim \Omega_{\frak L} =0$, 
and the $C^*$-algebra $\ARL$ is the $C^*$-algebra ${\Cal O}_{\frak L}$
associated with the $\lambda$-graph system $\frak L$.
\endroster

\heading 3. Condition (I) for  $C^*$-symbolic dynamical systems
\endheading

The notion of condition (I) for finite square matrices 
with entries in $\{0,1\}$ has been introduced in [CK]. 
The condition 
gives rise to the uniqueness of the associated Cuntz-Krieger algebras
under the canonical relations of the generating partial isometries.
The condition has been generalized by many authors to 
corresponding conditions for generalizations of the Cuntz-Krieger algebras,
 for instance,  infinite directed graphs ([KPRR]), infinite matrices with entries in $\{0,1\}$ ([EL]), Hilbert $C^*$-bimodules ([KPW]), etc. 
 (see also [Re], [Ka2],[Tom2],etc.).
 The condition (I) for $\lambda$-graph systems has been also defined in [Ma2]
 to prove the uniqueness of the $C^*$-algebra ${\Cal O}_{\frak L}$
 under the canonical relations of generators.
 In this section, 
 we will introduce the notion of condition (I) 
 for $C^*$-symbolic dynamical systems
 to prove 
 the uniqueness of the $C^*$-algebras
$\ARL$
under the relation $(\rho)$.
In [KPW], a condition called (I)-free has been introduced.
The condition is similar condition to our condition (I).
The discussions given in [KPW]
is also similar ones to ours in this section.
We will give 
complete descriptions in our discussions for the sake of completeness.
Throughout this paper, for a subset $F$ of a $C^*$-algebra $\B$,
we denote by $C^*(F)$ the $C^*$-subalgebra of $\B$ generated by $F$.

In what follows, $(\A,\rho,\Sigma)$ denotes a $C^*$-symbolic dynamical system
and $\Lambda$ the associated subshift $\Lambda_\rho$.
We denote by $\Lambda^k$ the set of admissible words $\mu$ of 
$\Lambda$ with length $| \mu | =k$.
Put
$
\Lambda^* = \cup_{k=0}^{\infty}\Lambda^k,
$
where $\Lambda^0 $ denotes the empty word. 
Let $S_\alpha, \alpha \in \Sigma$ 
be the partial isometries in $\ARL$ satisfying the relation $(\rho)$
 in Proposition 2.1.
For 
$\mu =(\mu_1,\dots,\mu_k) \in \Lambda^k$, 
we put
$S_\mu = S_{\mu_1}\cdots S_{\mu_k}$
and
$\rho_\mu = \rho_{\mu_k}\circ \cdots \circ \rho_{\mu_1}$.
In the algebra $\ARL$, 
we set
$$
\align
 {\Cal F}_{\rho} 
= &  C^*( S_{\mu}xS_{\nu}^{*}:
 \mu, \nu \in \Lambda^{*}, |\mu| = |\nu|, x \in \A), \\
 {\Cal F}_{\rho}^k 
= &  C^*( S_{\mu}xS_{\nu}^{*}: \mu, \nu \in \Lambda^{k}, x \in \A), 
          \text{ for } k \in \Zp \quad \text{ and}\\
{\Cal D}_{\rho} 
= &  C^*( S_\mu x S_{\mu}^*, \mu \in \Lambda^*, x \in \A). \\
\endalign
$$
The identity 
$
S_\mu x S_\nu^* 
= \sum_{\alpha \in \Sigma}S_{\mu \alpha} \rho_\alpha(x) S_{\nu \alpha}^*
$
for $x \in \A$ and $\mu, \nu \in \Lambda^k$ holds
so that the algebra
$
{\Cal F}_{\rho}^k
$ is embedded into the algebra
$
{\Cal F}_{\rho}^{k+1}
$
such that 
$
\cup_{k\in \Zp}{\Cal F}_{\rho}^k
$ is dense in ${\Cal F}_{\rho}$.
The gauge action $\hat{\rho}$ 
of the circle group ${\Bbb T} = \{ z \in {\Bbb C} \mid |z| = 1 \}$
on $\ARL$ is defined by
$\hat{\rho}_z(x) = x$ for $ x\in \A$
and $\hat{\rho}_z(S_\alpha) = z S_\alpha$
for $\alpha \in \Sigma$.
The fixed point algebra of $\ARL$ under $\hat{\rho}$ is denoted by
$(\ARL)^{\hat{\rho}}$. 
Let $\E_\rho: \ARL \longrightarrow (\ARL)^{\hat{\rho}}$ be the conditional expectaton defined by
$$
\E_\rho(X) = \int_{z \in \Bbb T} \hat{\rho}_z(X) dz, \qquad X \in \ARL.
$$
It is routine to check that 
$(\ARL)^{\hat{\rho}} = {\Cal F}_{\rho}.
$

Let $\B$ be a unital $C^*$-algebra.
Suppose that there exist an injective homomorphism
$\pi: \A \longrightarrow \B$ preserving their units
and 
a family $s_\alpha \in \B, \alpha \in \Sigma$ 
of partial isometries
satisfying 
$$
\sum_{\beta \in \Sigma} s_{\beta}s_{\beta}^* =1,\qquad
s_\alpha^* \pi(x) s_\alpha = \pi(\rho_\alpha(x)),
\qquad
\pi(x) s_\alpha s_\alpha^* =  s_\alpha s_\alpha^*  \pi(x)  
$$
for all $x \in \A$ and $\alpha \in \Sigma.$
 Put
 $\widetilde{\A} = \pi(\A) \subset $
 and
  $\tilde{\rho}_\alpha(\pi(x)) = \pi(\rho_\alpha(x)), x \in \A$.
  We then have
\proclaim{Lemma 3.1}
The triple
$(\widetilde{\A},\tilde{\rho},\Sigma)$ is a $C^*$-symbolic dynamical system
such that the presented subshift 
$\Lambda_{\tilde{\rho}}$
is the same as the one $\Lambda (=\Lambda_{\rho})$
presented by
$(\A,\rho,\Sigma)$.
\endproclaim
Let
${\Cal O}_{\pi,s}$
be the $C^*$-subalgebra of $\B$ generated by  
$\pi(x)$ and $s_\alpha$ for $ x \in \A, \alpha \in \Sigma.$
In the algebra ${\Cal O}_{\pi,s}$, we set
$$
\align
{\Cal F}_{\pi,s} 
= &  C^*( s_\mu \pi(x) s_{\nu}^*: 
\mu, \nu \in \Lambda^*,|\mu| = |\nu|, x \in \A), \\
{\Cal F}_{\pi,s}^k 
= &  C^*( s_\mu \pi(x) s_{\nu}^*:
 \mu, \nu  \in \Lambda^k, x \in \A) 
 \text{ for } k \in \Zp \quad \text{ and}\\
{\Cal D}_{\pi,s} 
= &  C^*( s_\mu \pi(x) s_{\mu}^*: \mu \in \Lambda^*, x \in \A).
\endalign
$$
By the universality of the algebra $\ARL$,
the correspondence
$$ 
x \in \A \longrightarrow \pi(x) \in \widetilde{A},
\qquad
S_\alpha \longrightarrow s_\alpha, \ \ \alpha \in \Sigma
$$
extends to a surjective homomorphism
$\tilde{\pi}:\ARL \longrightarrow {\Cal O}_{\pi,s}$. 
\proclaim{Lemma 3.2}
The restriction of $\tilde{\pi}$ to the subalgebra
${\Cal F}_{\rho}$
is an isomorphism from
${\Cal F}_{\rho}$ to ${\Cal F}_{\pi,s}$.
\endproclaim
\demo{Proof}
It suffices to show that $\tilde{\pi}$ is injective on
${\Cal F}_{\rho}^k$.
Suppose that 
$
\sum_{\mu,\nu\in \Lambda^k}s_\mu \pi(x_{\mu,\nu}) s_\nu^*=0
$
for
$
\sum_{\mu,\nu\in \Lambda^k}S_\mu x_{\mu,\nu} S_\nu^* 
\in
{\Cal F}_{\rho}
$
with
$x_{\mu,\nu}\in \A$.
For $\xi,\eta \in \Lambda^k$, it follows that 
$$
\pi(\rho_\xi(1)x_{\xi,\eta}\rho_\eta(1))
= s_\xi^* (\sum_{\mu,\nu\in \Lambda^k}s_\mu \pi(x_{\mu,\nu}) s_\nu^*)s_\eta
=0.
$$
As $\pi: \A \longrightarrow \B$ is injective,
one has
$\rho_\xi(1)x_{\xi,\eta}\rho_\eta(1)=0
$
so that 
$S_\xi x_{\xi,\eta}S_\eta^* =0$.
This implies that
$\sum_{\mu,\nu\in \Lambda^k}S_\mu x_{\mu,\nu} S_\nu^*=0.$
\qed
\enddemo

\noindent
{\bf Definition. }
A $C^*$-symbolic dynamical system $(\A,\rho,\Sigma)$
satisfies {\it condition }(I) if
there exists a unital increasing sequence
$$
\A_0 \subset \A_1 \subset \cdots \subset \A
$$
of $C^*$-subalgebras of $\A$ such that 
$\rho_\alpha(\A_l) \subset \A_{l+1}$ 
for all $l \in \Zp, \alpha \in \Sigma$
 and 
the union $\cup_{l\in \Zp} \A_l $ is dense in $\A$
and 
for $k,l\in \Bbb N$ with $k \le l$,
there exists a projection 
$
q_k^l \in {\Cal D}_{\rho}\cap {\A_l}^\prime 
(= \{ x \in {\Cal D}_\rho \mid x a = a x \ \text{ for } a \in \A_l \})
$
such that
\roster
\item"(i)" $q_k^l a \ne 0 $ for all nonzero $a \in \A_l$,
\item"(ii)" $q_k^l \phi_{\rho}^m(q_k^l) = 0 $ for all $m= 1,2,\dots, k,$
where 
$
 \phi_{\rho}^m(X) = 
 \sum_{\mu \in \Lambda^m} S_\mu X S_\mu^{*}.
$ 
\endroster
As the projection $q_k^l$ belongs to 
the diagonal subalgebra ${\Cal D}_\rho$ of ${\Cal F}_\rho$,
 the condition (I) of $(\A,\rho,\Sigma)$
 is intrinsically determined by $(\A,\rho,\Sigma)$
 by virtue of Lemma 3.2.
 
 If  a $\lambda$-graph system $\frak L$  over $\Sigma$
 satisfies condition (I),
 then $(\A_{\frak L}, \rho^{\frak L}, \Sigma)$ 
 satisfies condition (I) (cf. [Ma2;lemma 4.1]).

We now assume that $(\A,\rho,\Sigma)$ satisfies condition (I).
We set for $k \le l$
$$
 {\Cal F}_{\rho,l}^k 
 =   C^*( S_{\mu}xS_{\nu}^{*}: \mu, \nu \in \Lambda^{k}, x \in \A_l). 
$$
There exists an inclusion relation 
${\Cal F}_l^k \subset {\Cal F}_{l'}^{k'}$ for $k \le k'$ and $l \le l'$.
We put a projection
$Q_k^l = \phi_\rho^k(q_k^l)$  in ${\Cal D}_{\rho}.$

\proclaim{Lemma 3.3}
The map $
X \in  {\Cal F}_{\rho,l}^k
\longrightarrow  Q_k^l X Q_k^l \in Q_k^l {\Cal F}_{\rho,l}^k Q_k^l
$
is a surjective isomorphism.
\endproclaim
\demo{Proof}
As $q_k^l$ commutes with $\A_l$, 
for $x \in \A_l$ and $\mu,\nu \in \Lambda^k$, we have
$$
Q_k^l S_\mu x S_\nu^* 
= \sum_{\xi \in \Lambda^k} S_\xi q_k^l S_\xi^*S_\mu x S_\nu^*  
= S_\mu q_k^l S_\mu^*S_\mu x S_\nu^*
= S_\mu x q_k^l S_\nu^*,
$$
and similarly 
$S_\mu x S_\nu^* Q_k^l = S_\mu x q_k^l S_\nu^*$
so that $Q_k^l$ commutes with $S_\mu x S_\nu^*$.
Hence the map 
$
X \in  {\Cal F}_{\rho,l}^k
\longrightarrow  Q_k^l X Q_k^l \in Q_k^l {\Cal F}_{\rho,l}^k Q_k^l
$
defines a surjective homomorphism.
It remains to show that it is injective.
Suppose that 
$
Q_k^l (\sum_{\mu,\nu \in \Lambda^k}S_\mu x_{\mu,\nu} S_\nu^*)Q_k^l =0
$
for 
$
X = \sum_{\mu,\nu \in \Lambda^k}S_\mu x_{\mu,\nu} S_\nu^*
$
with
$
x_{\mu,\nu} \in \A_l.
$ 
For $\xi,\eta \in \Lambda^k$, one has 
$$
0 =
S_\xi S_\xi^* 
Q_k^l (\sum_{\mu,\nu \in \Lambda^k}S_\mu x_{\mu,\nu} S_\nu^*)Q_k^l 
S_\eta S_\eta^* 
=Q_k^l S_\xi x_{\xi,\eta} S_\eta^*,
$$
so that
$
0=
S_\xi^* Q_k^l S_\xi x_{\xi,\eta} S_\eta^* S_\eta
= S_\xi^* S_\xi  q_k^l \rho_\xi(1) x_{\xi,\eta} S_\eta^*S_\eta
= q_k^l\rho_\xi(1) x_{\xi,\eta} \rho_\eta(1).
$
Hence
$\rho_\xi(1) x_{\xi,\eta}\rho_{\eta}(1) =0$
by condition (I).
Thus  
$S_\xi x_{\xi,\eta} S_\eta^* =0$,
so that
$  \sum_{\xi,\eta\in \Lambda^k} S_\xi x_{\xi,\eta} S_\eta^* =0 $.
\qed
\enddemo
\proclaim{Lemma 3.4}
$Q_k^l S_\mu Q_k^l =0$ for $\mu \in \Lambda^*$ with $|\mu| \le k \le l$.
\endproclaim
\demo{Proof}
By condition (I), we have 
$Q_k^l \phi_\rho^m(Q_k^l) =0$ for $1 \le m \le k$.
For $\mu \in \Lambda^*$ with $|\mu | \le k$, one has 
$
\phi_\rho^{|\mu |}(Q_k^l)S_\mu = S_\mu Q_k^l S_\mu^* S_\mu = S_\mu Q_k^l.$
Hence we have
$
0 = 
Q_k^l \phi_\rho^{|\mu |}(Q_k^l)S_\mu = Q_k^l S_\mu Q_k^l.
$
\qed
\enddemo
As a result, we have
\proclaim{Lemma 3.5}
The projections $Q_k^l$ in ${\Cal D}_{\rho}$
satisfy the following conditions:
\roster
\item"(a)" $Q_k^l F - F Q_k^l $ converges to $0$ as $k,l \rightarrow \infty$
for $F \in 
{\Cal F}_{\rho}.
$
\item"(b)" $\| Q_k^l F \| $ converges to $ \| F \|$ as $k,l \rightarrow \infty$
for $F \in 
{\Cal F}_{\rho}.
$
\item"(c)"
$Q_k^l S_\mu Q_k^l =0$ for $\mu \in \Lambda^*$ with $|\mu| \le k \le l$.
\endroster
\endproclaim

We note that 
$ Q_k^l S_\mu Q_k^l =0$ if and only if $Q_k^l S_\mu Q_k^l S_\mu^*=0$.
Since 
$
Q_k^l S_\mu Q_k^l S_\mu^*
$ belongs to the algebra
$
{\Cal F}_{\rho},
$
the condition
$
Q_k^l S_\mu Q_k^l  =0
$
is determined in the algebraic structure of 
$
{\Cal F}_{\rho}.
$
As the restriction of 
$\tilde{\pi}: \ARL \longrightarrow {\Cal O}_{\pi,s}$
to
${\Cal F}_\rho$
yields an isomorphism onto
${\Cal F}_{\pi,s}$,
by putting
$\widetilde{Q}_k^l = \tilde{\pi} (Q_k^l)$ 
we have

\proclaim{Lemma 3.6}
The  projections $\widetilde{Q}_k^l$
in ${\Cal D}_{\pi,s}$  satisfy the following conditions:
\roster
\item"(a')" $\widetilde{Q}_k^l F - F \widetilde{Q}_k^l $ converges to $0$ as $k,l \rightarrow \infty$
for $F \in 
{\Cal F}_{\pi,s}.
$
\item"(b')" $\| \widetilde{Q}_k^l F \| $ converges to $ \| F \|$ as $k,l \rightarrow \infty$
for $F \in 
{\Cal F}_{\pi,s}.
$
\item"(c')"
$\widetilde{Q}_k^l s_\mu \widetilde{Q}_k^l =0$ for $\mu \in \Lambda^*$ with $|\mu| \le k \le l$.
\endroster
\endproclaim

\proclaim{Proposition 3.7}
There exists a conditional expectation
$\E_{\pi,s}:{\Cal O}_{\pi,s} \longrightarrow {\Cal F}_{\pi,s}
$
such that
$\E_{\pi,s} \circ \tilde{\pi} =  \tilde{\pi}\circ \E_\rho.
$
\endproclaim
\demo{Proof}
Let ${\Cal P}_{\pi,s}$
be the $*$-subalgebra of ${\Cal O}_{\pi,s}$ 
generated algebraically by 
$\pi(x), s_\alpha$ for $x \in \A, \alpha \in \Sigma$.
Then any  
$ X \in {\Cal P}_{\pi,s}$ can be written as a finite sum
$$
X = \sum_{|\nu| \ge 1} X_{-\nu} s_\nu^* + X_0 +  
    \sum_{|\mu| \ge 1} s_\mu X_{\mu}
 \quad
 \text{ for some }   X_{-\nu}, X_0, X_\mu \in {\Cal F}_{\pi,s}. 
$$
Thanks to the previous lemma and a usual argument of [CK],
 the element 
$X_0 \in {\Cal F}_{\pi,s}$ 
is unique for $X \in {\Cal P}_{\pi,s}$
and the inequality
$ 
\| X_0 \| \le \| X \|
$
holds.   
The map 
$X \in {\Cal P}_{\pi,s} \rightarrow X_0 \in {\Cal F}_{\pi,s}$
can be extended to the desired expectation
$\E_{\pi,s}:{\Cal O}_{\pi,s} \longrightarrow {\Cal F}_{\pi,s}.
$
\qed
\enddemo

Therefore we have
\proclaim{Theorem 3.8}
Assume that $(\A,\rho,\Sigma)$ 
satisfies condition (I).
The homomorphism $\tilde{\pi}:\ARL \longrightarrow {\Cal O}_{\pi,s}$
defined by
$$
\tilde{\pi}(x) = \pi(x),\quad x \in \A, \qquad
\tilde{\pi}(S_\alpha) = s_\alpha, \quad \alpha \in \Sigma.
$$
becomes a surjective isomorphism, and hence
the $C^*$-algebras $\ARL$ and ${\Cal O}_{\pi,s}$ 
are canonically isomorphic through $\tilde{\pi}$.
\endproclaim
\demo{Proof}
The map $\tilde{\pi}: {\Cal F}_{\rho} \rightarrow {\Cal F}_{\pi,s}$
is isomorphic and satisfies 
$\E_{\pi,s} \circ \tilde{\pi} =  \tilde{\pi}\circ \E_\rho.
$
Since 
$\E_\rho:\ARL \longrightarrow {\Cal F}_{\rho}$
is faithful,
a routine argument shows that the homomorphism 
$\tilde{\pi}:\ARL \longrightarrow {\Cal O}_{\pi,s}$
is actually an isomorphism.
\qed
\enddemo
Hence the following uniqueness of the $C^*$-algebra $\ARL$ holds.
\proclaim{Theorem 3.9}
Assume that $(\A,\rho,\Sigma)$ 
satisfies condition (I).
The $C^*$-algebra $\ARL$ is the unique $C^*$-algebra subject to the relation 
$(\rho)$.
This means that if there exist a unital $C^*$-algebra $\B$ and an
injective homomorphism $\pi: \A \longrightarrow \B$ and a family 
$s_\alpha \in \B, \alpha \in \Sigma$
of nonzero partial isometries satisfying the folloing relations:
$$
\sum_{\beta \in \Sigma} s_{\beta}s_{\beta}^* =1,\qquad
s_\alpha^* \pi(x) s_\alpha = \pi(\rho_\alpha(x)),
\qquad
\pi(x )s_\alpha s_\alpha^* =  s_\alpha s_\alpha^* \pi( x ) 
$$
for all $x \in \A$ and $\alpha \in \Sigma$,
then the correspondence 
$$
 x \in \A \longrightarrow \pi(x) \in \B, \qquad
S_\alpha \longrightarrow s_\alpha \in \B
$$ 
extends to an  isomorphism
$\tilde{\pi}$ from $ \ARL$ onto the $C^*$-subalgebra
${\Cal O}_{\pi,s}$ of $\B$ 
generated by $\pi(x), x \in \A$ and $s_\alpha, \alpha \in \Sigma$.
\endproclaim
As a corollary we have
\proclaim{Corollary 3.10}
Assume that $(\A,\rho,\Sigma)$ 
satisfies condition (I).
For any nontrivial ideal ${\Cal I} $
of $\ARL$, one has 
${\Cal I} \cap \A \ne \{ 0\}.$  
\endproclaim
\demo{Proof}
Suppose that  
${\Cal I} \cap \A = \{ 0 \}$.
Hence $S_\alpha \not\in {\Cal I}$ 
for all $\alpha \in \Sigma$.
By Theorem 3.9, the quotient map
$q: \ARL \longrightarrow \ARL/{\Cal I}$
must be injective so that ${\Cal I}$ is trivial.
\qed
\enddemo

Let
$\lambda_\rho: \A \rightarrow \A$
be the completely positive map on $\A$ defined by
$\lambda_\rho(x) = \sum_{\alpha \in \Sigma}\rho_\alpha(x)$ 
for $x \in \A$.

\noindent
{\bf Definition.}
$(\A,\rho,\Sigma)$
is said to be  {\it irreduible } if
there exists no nontrivial ideal of $\A$ invariant under $\lambda_\rho$.

Therefore we have
\proclaim{Corollary 3.11}
Assume that $(\A,\rho,\Sigma)$ 
satisfies condition (I).
If $(\A,\rho,\Sigma)$ is irreducible,
the $C^*$-algebra $\ARL$ is simple.
\endproclaim

\heading 4. Quotients of  $C^*$-symbolic dynamical systems
\endheading
In this section, 
we will study ideal structure of the $C^*$-symbolic crossed products
$\ARL$, related to quotients of $C^*$-symbolic dynamical systems.
The ideal structure of $C^*$-algebras of Hilbert $C^*$-bimodules 
has been studied in Kajiwara, Pinzari and Watatani's paper [KPW] 
(cf. [Kat3]). 
Their paper is written in the language of Hilbert $C^*$-bimodules.
In this section we will directly study ideal structure of 
the $C^*$-symbolic crossed products
$\ARL$ by using the language of $C^*$-symbolic dynamical systems.
We fix a $C^*$-symbolic dynamical system $(\A,\rho,\Sigma)$.

An ideal $J$ of $\A$ is said to be $\rho$-{\it invariant } if
$\rho_\alpha (J) \subset J$ for all $\alpha \in \Sigma$.
It is said to be {\it saturated } if $\rho_\alpha(x) \in J$ for all $\alpha \in \Sigma$ implies $x \in J$.
\proclaim{Lemma 4.1} Let $J$ be an ideal of $\A$.
\roster
\item"(i)" $J$ is $\rho$-invariant if and only if $\lambda_\rho(J) \subset J$.
\item"(ii)" $J$ is saturated if and only if 
$\lambda_\rho(a)\in J $ for  $ 0 \le a\in \A$ implies $a \in J$.
\endroster
\endproclaim
\demo{Proof}
(i) Suppose that $J$ satisfies $\lambda_\rho(J) \subset J$.
For $x \in J$ one has 
$
\lambda_\rho(x^*x) 
\ge \rho_\alpha(x^*x) =\rho_\alpha(x)^*\rho_\alpha(x)
$
 so that 
$\rho_\alpha(x)^*\rho_\alpha(x) \in J$
because ideal is hereditary.
Hence $\rho_\alpha(x)$ belongs to $J$.
The only if part is clear.

(ii) Suppose that $J$ is saturated and $\lambda_\rho(a)\in J $ for 
$ 0 \le a\in \A$. 
Since $\lambda_\rho(a) \ge \rho_\alpha(a)$ and $J$ is hereditary,
one has $a \in J$.
Conversely suppose that $x\in \A$ satisfies 
$\rho_\alpha(x) \in J$ for all $\alpha \in \Sigma$.
As 
$
\lambda_\rho(x^*x) 
= \sum_{\alpha\in \Sigma}\rho_\alpha(x)^* \rho_\alpha(x),
$
$
\lambda_\rho(x^*x) 
$ belongs to $J$.
Hence the condition of the if part implies that $x^*x \in J$
so that $x \in J$.
\qed
\enddemo
Let $J$ be a $\rho$-invariant saturated ideal of $\A$.
We denote by ${\Cal I}_J$  the ideal of $\ARL$ generated by $J$.
\proclaim{Lemma 4.2}
The ideal  ${\Cal I}_J$
is the closure of linear combinations of elements of the form
$
 S_\mu c_{\mu,\nu} S_\nu^*
$
for
$
 c_{\mu,\nu} \in J.
$
\endproclaim
\demo{Proof}
Elements $x$ and $ y $ of $\ARL$ 
are approximated by finite sums of elements of the form
$
 S_\mu a_{\mu,\nu} S_\nu^*
$
and
$
 S_\xi b_{\xi,\eta} S_\eta^*
$
for
 $ a_{\mu,\nu}, b_{\xi,\eta} \in \A
$
respectively.
Hence $x c y$ is approximated by elemnts of the form
$$
\sum_{\mu, \nu} S_\mu a_{\mu,\nu} S_\nu^* \cdot c \cdot
\sum_{\xi, \eta} S_\xi b_{\xi,\eta} S_\eta^*
=
\sum_{\mu, \nu, \xi, \eta} S_\mu a_{\mu,\nu} S_\nu^* c
 S_\xi b_{\xi,\eta} S_\eta^*.
$$
In case of $|\nu| \ge |\xi|$,
one has $\nu = \bar{\nu} \nu'$ with $|\bar{\nu}| = | \xi |$
so that
$$
S_\nu^* c S_\xi
=
\cases
S_{\nu'}^* \rho_{\bar{\nu}}(c) & \text{ if } \bar{\nu} = \xi,\\ 
0 & \text{ otherwise. }
\endcases
$$
Hence 
$S_\mu a_{\mu,\nu} S_\nu^* c
 S_\xi b_{\xi,\eta} S_\eta^*
 $
is 
$
S_\mu a_{\mu,\nu} \rho_{\nu'}( \rho_{\bar{\nu}}(c) b_{\xi,\eta}) 
S_{\eta \nu'}^*
$
or zero.
Since $J$ is $\rho$-invariant,
it is of the form 
$S_\mu c_{\mu,\nu} S_{\eta \nu'}^*$ for some $c_{\mu,\nu}\in J$
or zero.
The argument in case of $|\nu| \le |\xi|$ is similar.
Since the ideal 
${\Cal I}_J$
is the closure of  elements of the form
$
\sum_{i=1}^n x_i c_i y_i
$
for
$ x_i, y_i \in \ARL
$
and
$ c_i \in J,
$
the assertion is proved.
\qed
\enddemo
We set 
$$
\align
{\Cal D}_J 
&= C^*(S_\mu c_\mu S_\mu^* : \mu \in \Lambda^*, c_\mu \in J),\\
{\Cal D}_J^k 
&= C^*(S_\mu c_\mu S_\mu^* : \mu \in \Lambda^*,|\mu| \le k, c_\mu \in J)
\quad \text{ for } k \in \Zp.\\
\endalign
$$
\proclaim{Lemma 4.3}
\roster
\item"(i)" $ {\Cal D}_J  = {\Cal I}_J \cap {\Cal D}_\rho$ and hence
${\Cal D}_J \cap \A =  {\Cal I}_J \cap \A.$
\item"(ii)" ${\Cal D}_J^k \cap \A = J$ for $k\in \Zp.$
\endroster
\endproclaim
\demo{Proof}
(i)
Since the elements of the finite sum 
$\sum_{\mu} S_\mu c_\mu S_\mu^*$ for $c_\mu \in J$ 
are contained in ${\Cal I}_J \cap {\Cal D}_\rho$,
the inclusion relation
$ {\Cal D}_J  \subset  {\Cal I}_J \cap {\Cal D}_\rho$ 
is clear.
Let
${\Cal I}_J^{\text{alg}}$ and ${\Cal D}_J^{\text{alg}}$
be the algebraic linear spans of
$ S_\mu c_{\mu,\nu} S_\nu^*$ for $c_{\mu,\nu} \in J$
and
$ S_\mu c_{\mu} S_\mu^*$ for $c_{\mu} \in J$
respectively. 
For any $x\in {\Cal I}_J \cap {\Cal D}_\rho$
take $x_n \in {\Cal I}_J^{\text{alg}}$ 
such that 
$\| x_n - x \| \rightarrow 0$.
Let
$\E_\rho :\ARL \rightarrow {\Cal F}_\rho$ 
be the conditional expectation defined previously,
and $\E_D: {\Cal F}_\rho \rightarrow {\Cal D}_\rho$
the conditional expectation defined by taking diagonal elements.
The composition 
$\E_{{\Cal D}_\rho} = \E_D \circ \E_\rho $ 
is the conditional expectation from $\ARL$ to ${\Cal D}_\rho$
that satisfies 
$\E_{{\Cal D}_\rho}({\Cal I}_J^{\text{alg}}) = {\Cal D}_J^{\text{alg}}$.
Since
$\E_{{\Cal D}_\rho}(x) = x$ and 
the inequality 
$
\| x - \E_{{\Cal D}_\rho}(x_n) \| \le \| x- x_n \|
$
holds,
 $x$ belongs to the closure ${\Cal D}_J$
of ${\Cal D}_J^{\text{alg}}$.
Hence we have
$ {\Cal I}_J \cap {\Cal D}_\rho \subset {\Cal D}_J$
so that
${\Cal D}_J = {\Cal I}_J \cap {\Cal D}_\rho$.
As $\A$ is a subalgebra of ${\Cal D}_\rho$,
the equality
${\Cal D}_J \cap \A =  {\Cal I}_J \cap \A$
holds.

(ii)
An element 
$x \in {\Cal D}_J^k$ is of the form
$\sum_{|\mu | \le k} S_\mu c_\mu S_\mu^*$ for $c_\mu \in J$.
As 
$
S_\nu c_\nu S_\nu^* 
= \sum_{\alpha \in \Sigma}S_{\nu \alpha}\rho_\alpha(c_\nu)S_{\nu \alpha}^*
$
and $J$ is $\rho$-invariant,
$x$ can be written as 
$
x= \sum_{|\nu | = k} S_\nu c_\nu S_\nu^* 
$ for $c_\nu \in J$,
and the element
$
\lambda_\rho^k(x) = \sum_{|\nu | = k} \rho_\nu(1) c_\nu \rho_\nu(1)
$
belongs to $J$.
Further suppose that $x$ is an element of $\A$.
Since $J$ is saturated, by Lemma 4.1, one has
$x \in J$.
Hence 
the inclusion relation
$\A \cap {\Cal D}_J^k \subset J$ holds.
The converse inclusion relation is clear so that
$\A \cap {\Cal D}_J^k = J$.   
\qed
\enddemo
\proclaim{Lemm 4.4}
$\A \cap {\Cal D}_J = J$.
\endproclaim
\demo{Proof}
Since the inclusion relation
$\A \cap {\Cal D}_J \supset J$ is clear,
there exists a natural surjective homomorphism
from
$\A / J$ onto $\A / \A \cap {\Cal D}_J$.
For  an element $a$ of a $C^*$-algebra $\B$, 
we denote by $\| [a]_{\B /I} \|$ 
the norm of the quotient image 
$[a]_{\B /I}$ of $a$ in the quotient
$\B / I$ of $\B$ by an ideal $I$.
As the inclusion $ \A \hookrightarrow {\Cal D}_\rho$
induces the  inclusions both
$\A / \A \cap {\Cal D}_J \hookrightarrow {\Cal D}_\rho / {\Cal D}_J$
and
$\A / \A \cap {\Cal D}_J^k \hookrightarrow {\Cal D}_\rho / {\Cal D}_J^k$,
one has for $a \in \A$
$$
\| [a]_{\A / \A \cap {\Cal D}_J } \| 
= \| [a]_{{\Cal D}_\rho / {\Cal D}_J}\|,
\qquad
\| [a]_{\A / \A \cap {\Cal D}_J^k } \| 
= \| [a]_{{\Cal D}_\rho / {\Cal D}_J^k}\|.
$$
Note that ${\Cal D}_J$ is the inductive limit of 
${\Cal D}_J^k, k=0,1,\dots.$
It then follows that
$$
\| [a]_{{\Cal D}_\rho / {\Cal D}_J}\|
 = \dist(a, {\Cal D}_J) 
 = \lim_{k \to \infty} \dist(a, {\Cal D}_J^k) 
 = \lim_{k \to \infty} \| [a ]_{{\Cal D}_\rho / {\Cal D}_J^k} \| 
 = \lim_{k \to \infty} \| [a]_{\A / \A \cap {\Cal D}_J^k } \| 
 $$
and hence  
$
\| [a]_{\A / \A \cap {\Cal D}_J } \| = \| [a]_{\A / J } \|
$
by Lemma 4.3 (ii).
Thus the quotient map
$ \A / J \rightarrow \A / \A \cap {\Cal D}_J$
is isometric  so that 
$\A \cap {\Cal D}_J = J$.
\qed
\enddemo
By Lemma 4.3 and Lemma 4.4, one has
\proclaim{Proposition 4.5}
$
 {\Cal I}_J \cap \A = J.
$
 \endproclaim
We will now consider quotient $C^*$-symbolic dynamical systems.
Let $J$ be a $\rho$-invariant saturated ideal of $\A$.
We set 
$
\Sigma_J = \{ \alpha \in \Sigma \mid \rho_\alpha(1) \not\in J \}.
$
We denote by $[x]$ the class of $x \in \A$ in the quotient $\A / J$.
Put 
$$
\rho_\alpha^J([x]) = [ \rho_\alpha(x) ] \qquad \text{ for } [x] \in \A / J,
\quad \alpha \in \Sigma_J.
$$
As $J$ is $\rho$-invariant and saturated,
$\rho_\alpha^J$ is well-defined
and 
the family 
$ 
\{ \rho_\alpha^J \}_{\alpha \in \Sigma_J}
$
is a faithful and essential endomorphisms of $\A / J$.
We call the $C^*$-symbolic dynamical system $(\A/J, \rho^J,\Sigma_J)$
the quotient of $(\A, \rho, \Sigma)$ by the ideal $J$.
We denote by $\Lambda_{J}$ the associated subshift for the quotient
$(\A/J, \rho^J,\Sigma_J)$.

\noindent
{\bf Definition.}
A $C^*$-symbolic dynamical system 
$(\A,\rho,\Sigma)$ is said to satisfy {\it condition } (II)
if for any proper $\rho$-invariant saturated ideal $J$ of $\A$,
the quotient $C^*$-symbolic dynamical system $(\A/J, \rho^J,\Sigma_J)$
satisfies condition (I).

Let $\I$ be a proper ideal of $\ARL$.
Put $J_\I := \I \cap \A$.
\proclaim{Lemma 4.6}
\roster
\item"(i)" If $(\A,\rho, \Sigma)$ satisfies condition (I), 
then
$J_\I$ 
is a proper $\rho$-invariant saturated ideal of $\A$.
We then have $J_{\I_J} = J$.
\item"(ii)" If $(\A,\rho, \Sigma)$ satisfies condition (II), 
then the $C^*$-symbolic crossed product
$
(\A / {J_\I}) \rtimes_{\rho^{J_\I}} \Lambda_{J_\I}
$
is canonically isomorphic to the quotient algebra
$\ARL / \I$. 
\endroster
\endproclaim
\demo{Proof}
(i) By condition (I), $J_\I$ is a nozero ideal of $\A$,
that is $\rho$-invariant.
If $\rho_\alpha(x) $ belongs to $J_\I$ for all $\alpha \in \Sigma$, 
the identity 
$x = \sum_{\alpha \in \Sigma} S_\alpha \rho_\alpha(x) S_\alpha^*$
implies $x \in \I$, so that $J_\I$ is saturated.
The equality $J_{\I_J} = J$ follows from Proposition 4.5.

(ii)
Let $\pi_\I: \ARL \rightarrow \ARL / \I$
be the quotient map.
Put $s_\alpha = \pi_\I(S_\alpha)$.
Then $\alpha \in \Sigma_{J_\I}$ if and only if $s_\alpha \ne 0$.
The following relations
$$
\sum_{\beta \in \Sigma_{J_\I}}s_\beta s_\beta^* = 1,
\qquad s_\alpha^*\pi_\I(x) s_\alpha = \pi_\I(\rho_\alpha(x)) 
\qquad
\pi_\I(x) s_\alpha s_\alpha^* =  s_\alpha s_\alpha^* \pi_\I(x)
$$
for $x \in \A, \alpha \in \Sigma_{J_\I}$
hold.
As 
$(\A/J_\I, \rho^{J_\I},\Sigma_{J_\I})$
satisfies condition (I),
the uniqueness of the $C^*$-symbolic crossed product
$
(\A / {J_\I}) \rtimes_{\rho^{J_\I}} \Lambda_{J_\I}
$
yields a canonical isomorphism to the quotient algebra
$\ARL / \I$.
\qed
\enddemo
Let $\I_{J_\I}$ be the ideal of $\ARL$ generated by $J_\I$.
Since $J_\I \subset \I$, the inclusion relation
$\I_{J_\I} \subset \I$ is clear.
\proclaim{Lemma 4.7}
 If $(\A,\rho, \Sigma)$ satisfies condition (II), 
then there exists a canonical isomorphism from 
$
(\A / {J_\I}) \rtimes_{\rho^{J_\I}} \Lambda_{J_\I}
$
to  the quotient algebra
$\ARL / {\I_{J_\I}}$. 
\endproclaim
\demo{Proof}
Take an arbitrary element $x \in \A$.
If $x \in J_\I$, then $x \in \I_{J_\I}$.
Conversely $x \in \I_{J_\I}$ implies $x \in J_\I$ by Proposition 4.5.
Hence 
$x \in J_\I$ if and only if $x \in \I_{J_\I}$.
For $\alpha \in \Sigma$,
we have
$S_\alpha \in \I_{J_\I}$ 
if and only if 
$S_\alpha^*S_\alpha \in \I_{J_\I}\cap \A$.
By Proposition 4.5, the latter condition is equivalent to the condition 
$\rho_\alpha(1) \in J_\I$.
We know that 
$\alpha \not\in \Sigma_{J_\I}$ 
if and only if $S_\alpha \in \I_{J_\I}$.
By the uniqueness of the algebra
$
(\A / {J_\I}) \rtimes_{\rho^{J_\I}} \Lambda_{J_\I},
$
it is canonically isomorphic to the quotient algebra
$\ARL / {\I_{J_\I}}$. 
  \qed
\enddemo
\proclaim{Proposition 4.8}
 Suppose that $(\A,\rho, \Sigma)$ satisfies condition (II). 
For a proper ideal $\I$ of $\ARL$,
let $\I_{J_\I}$ be the ideal of $\ARL$ generated by
$J_\I$.
Then we have
$
\I_{J_\I} = \I.
$
\endproclaim
\demo{Proof}
Since $\I_{J_\I} \subset \I$,
there exists a quotient map
$q_\I: \ARL / \I_{J_\I} \rightarrow \ARL / \I$.
By Lemma 4.6, and Lemma 4.7, there exist canonical isomorphisms
$$
\pi_1: (\A / {J_\I}) \rtimes_{\rho^{J_\I}} \Lambda_{J_\I}
 \rightarrow 
 \ARL / \I,
\qquad
 \pi_2: (\A / {J_\I}) \rtimes_{\rho^{J_\I}} \Lambda_{J_\I}
\rightarrow 
\ARL / \I_{J_\I}.
$$
Since $q_\I = \pi_1 \circ \pi_2^{-1}$,
it is ismorphism so that we have $\I_{J_\I} = \I$.
\qed
\enddemo

Therefore we have
\proclaim{Theorem 4.9}
 Suppose that $(\A,\rho, \Sigma)$ satisfies condition (II). 
There exists an inclusion preserving bijective correspondence between 
$\rho$-invariant saturated ideals of $\A$ and ideals of $\ARL$,
through the correspondences:
$J \rightarrow \I_J$ and $J_\I \leftarrow \I$.
\endproclaim

\heading 5. Pure infiniteness
\endheading
In this section we will show that 
the $C^*$-symbolic crossed product  $\ARL$
is purely infinite 
if $(\A,\rho,\Sigma)$ satisfies some conditions.

\noindent
{\bf Definition.}
A $C^*$-symbolic dynamical system
$(\A,\rho,\Sigma)$ is said to be {\it central} if the projections
$\{\rho_\mu(1) \mid \mu \in \Lambda^*\}$ contained 
in the center $Z_\A$ of $\A$.
It is said to be {\it commutative}
if $\A$ is commutative. 
Hence if $(\A,\rho,\Sigma)$ is central,
the inequality $\sum_{\alpha \in \Sigma}\rho_\alpha(1) \ge 1$ holds.
Let $\A_\rho$ be the $C^*$-subalgebra of $\A$ generated by the projections
$\rho_\mu(1), \mu \in \Lambda^*$.

\proclaim{Lemma 5.1}
Assume that  
$(\A,\rho,\Sigma)$ is central.
Then 
there exists a $\lambda$-graph system
${\frak L}_{\rho}$
over $\Sigma$ such that 
the presented subshift
$\Lambda_{{\frak L}_{\rho}}$
coincides with the subshift $\Lambda$
 presented by $(\A,\rho,\Sigma)$,
and there exists a unital embedding of
${\Cal O}_{{\frak L}_{\rho}}$
into $\ARL$.
\endproclaim
\demo{Proof}
Put 
$
\A_{\rho,0} = {\Bbb C}.
$
For $l \in \Zp$, we define the $C^*$-algebra 
$
\A_{\rho,l+1} 
$
to be  the $C^*$-subalgebra of $\A$ 
generated by the elements 
$\rho_\alpha (x)$ for $\alpha \in \Sigma, x \in \A_{\rho,l}$.
Hence the $C^*$-algebra $\A_\rho$ is generated by
$\cup_{l = 0}^{\infty}\A_{\rho,l}$.
Then  
$({\A_\rho},\rho,\Sigma)$ is a $C^*$-symbolic dynamical system
such that 
$\A_\rho$ is commutative and AF,
so that 
there exists a $\lambda$-graph system
${\frak L}_{\rho}$
over $\Sigma$ such that 
$\A_\rho = \A_{{\frak L}_{\rho}}$.
The presented subshift
$\Lambda_{{\frak L}_\rho}$
coincides with the subshift $\Lambda$.
It is easy to see that there exists a unital embedding of
${\Cal O}_{{\frak L}_{\rho}}$
into $\ARL$ by their universalities.
\qed
\enddemo

In the rest of  this section we assume that $(\A,\rho,\Sigma)$ satisfies condition (I).

\noindent
{\bf Definition.}
$(\A,\rho,\Sigma)$ is said to be {\it effective}
if for $l \in \Zp$ and a nonzero positive element $a \in \A_l$,
there exist $K \in \Bbb N$ and a nonzero positive element $b \in \A_\rho$
such that
$$
\sum_{\mu \in \Lambda^K} \rho_\mu(a) \ge b \tag 5.1
$$  
where $\A_l$ is a $C^*$-subalgebra of $\A$ appearing in the definition of condition (I).

In what follows, we assume that  $(\A,\rho,\Sigma)$ is effective, 
and central.
 Let ${\frak L} ={\frak L}_{\rho}$ be the $\lambda$-graph system associated to $(\A,\rho,\Sigma)$ as in Lemma 5.1.
 We further assume that the algebra 
 ${\Cal O}_{\frak L}$ is simple, purely infinite.
 In [Ma2], [Ma3], conditions that the algebra ${\Cal O}_{\frak L}$ becomes simple, purely infinite is studied.
   
\proclaim{Lemma 5.2}
For $k\le l \in \Zp$ 
and a nonzero positive element 
$a \in 
{\Cal F}_{\rho,l}^k
$,
there exists an element $V  \in \ARL$
such that
$
V a V^* =1.
$ 
\endproclaim
\demo{Proof}
An element 
$
a \in {\Cal F}_{\rho,l}^k
$ is of the form
$ a = \sum_{\mu,\nu \in \Lambda^k} S_\mu a_{\mu,\nu}  S_\nu^*$ 
for some $ a_{\mu,\nu} \in \A_l$ such that
$S_\mu^*  a S_\nu = a_{\mu,\nu}$.
Since $a$ is a nonzero positive element, 
there exists $\xi \in \Lambda^k$
such that
$S_\xi^* a S_\xi (= a_{\xi,\xi}) \ne 0.$
As we are assuming that 
$(\A,\rho, \Sigma)$ is effective,
there exists $K \in \Bbb N$ and a nonzero positive element $b \in \A_\rho$
$$
\sum_{\mu \in \Lambda^K} \rho_\mu(S_\xi^* a S_\xi ) \ge b. 
$$  
Put $T = \sum_{\mu \in \Lambda^K} S_\mu \in \ARL.$
One has $T^* S_\xi^* a S_\xi T \ge b$.
Now $b \in \A_\rho \subset {\Cal O}_{\frak L}$ 
and ${\Cal O}_{\frak L}$ is simple, purely infinite.
We may find $V_0 \in {\Cal O}_{\frak L}$ 
such that $V_0 b V_o^* = 1$
so that 
$V_0 T^* S_\xi^* a S_\xi T V_0^* \ge 1$.
Hence there exists $V \in \ARL$ such that 
$V a V^* =1.$
\qed
\enddemo

\proclaim{Lemma 5.3}
Keep the above situation.
We may take $V \in \ARL$ in the preceding lemma such as
$
V a V^* =1
$
and
$
 \| V \| < \| a \| ^{-\frac{1}{2}} + \epsilon
$ 
for a given 
$\epsilon >0$.
\endproclaim
\demo{Proof}
We may assume that $\| a \| =1$ 
and there exists $p\in Sp(a)$ such that $0 < p < 1$.
Take $ 0 < \epsilon < \frac{1}{2}$
such that  $\epsilon < 1-p$.
Define a function $f \in C([0,1])$ 
by setting
$$
f(t) =
\cases
0 & (0\le t \le 1-\epsilon) \\
1 - \epsilon^{-1}(1-t) & 1-\epsilon < t \le 1 \\
\endcases
$$
Put $b= f(a)$, that is not invertible.
By Lemma 5.2, there exists $V \in \ARL$ 
such that
$V b V^* =1$.
We set $S = b^{\frac{1}{2}} V^* $
and
$P= S S^*$.
Then $S$ is a proper isometry
such that $P \le \| V \| b$.
As $P \le E_a([1-\epsilon,1])$,
$E_a([1-\epsilon,1])$ is the spectral measure of $a$ 
for the interval $[1-\epsilon,1]$,
one has
$P a P \ge (1-\epsilon) P$. 
Put
$D = S^* a S$ so that
$
D \ge S^* (1-\epsilon) PS = (1 -\epsilon)1.
$
Hence $D$ is invertible. 
Set
$V_1 = D^{-\frac{1}{2}} S^*$.
Then one sees that 
$ 
V_1 a V_1^* = 1
$ 
and
$
\| V_1 \| < (1 - \epsilon )^{-\frac{1}{2}} < 1 + \epsilon. 
$
\qed
\enddemo
Let $\E_\rho: \ARL \rightarrow {\Cal F}_{\rho}$ be 
the conditional expectation defined in Section 3.
\proclaim{Lemma 5.4}
For a nonzero $X \in \ARL$ and $\epsilon >0$,
there exists a projection $Q \in {\Cal D}_\rho$
and a nonzero positive element $Z \in {\Cal F}_{\rho,k}^l$ for some $k \le l$
such that
$$
\| Q X^*X Q -Z\| < \epsilon, \qquad 
\| \E_\rho (X^*X)\| -\epsilon < \| Z \| < \| \E_\rho (X^*X) \| + \epsilon.
$$
\endproclaim
\demo{Proof}
We may assume that 
$\| \E_\rho (X^*X) \| =1.$
Let
${\Cal P}_{\rho}$ be the 
$*$-algebra generated algebraically by 
$S_{\alpha},\alpha \in \Sigma$ and 
$x \in \A.$
For any $0 < \epsilon < \frac{1}{4},$
find 
$0\le Y \in 
{\Cal P}_{\rho}$ 
such that
$
\| X^*X - Y \| < \frac{\epsilon}{2}
$
so that 
$
\| \E_\rho(Y) \| >  1 - \frac{\epsilon}{2}.
$
As in the discussion in [Ma3;Section 3], the element
$Y$ is expressed 
as 
$$
Y = \sum_{|\nu | \ge 1}Y_{-\nu} S_{\nu}^* + Y_0 + 
    \sum_{|\mu | \ge 1}S_{\mu} Y_{\mu}
\qquad
\text{ for some }
\quad
Y_{-\nu}, Y_0, Y_{\mu}\in {\Cal F}_{\rho} \cap {\Cal P}_{\rho}.
$$
Take $k\le l$ large enough such that 
$
Y_{-\nu}, Y_0, Y_{\mu}\in {\Cal F}_{\rho,k}^l
$
for all
$
\mu, \nu
$ 
in the above expression.
Now 
$(\A,\rho,\Sigma)$ satisfies condition (I).
Take a sequence $Q_k^l \in {\Cal D}_\rho$
of projections as in Section 3.
As $\E_\rho(Y) = Y_0$ 
and $Q_k^l$ commutes with ${\Cal F}_{\rho,k}^l$,
it follows that by Lemma 3.5 (c), 
$
Q_k^l Y Q_k^l 
= Q_k^l \E_\rho(Y) Q_k^l.
$
Since $Q_k^l \E_\rho(Y) Q_k^l \in {\Cal F}_{\rho},$
there exists $0 \le Z \in {\Cal F}_{\rho,k'}^{l'}$ 
for some $k' \le l'$
such that 
$
\| Q_k^l \E_\rho(Y) Q_k^l -Z \| < \frac{\epsilon}{2}.
$
By Lemma 3.3,
we note 
$\| Q_k^l \E_\rho(Y) Q_k^l \| = \| \E_\rho(Y) \|$ so that 
$$
\| Z \| \ge \| \E_\rho(Y) \| - \frac{\epsilon}{2} > 1 - \epsilon
$$
and
$$
\| Z \| < \| Q_k^l \E_\rho(Y) Q_k^l \| + \frac{\epsilon}{2}
        \le \|  \E_\rho(X^* X) \| + \frac{\epsilon}{2} +\frac{\epsilon}{2}
        < 1 + \epsilon.
$$
\qed
\enddemo
Therefore we have
\proclaim{Theorem 5.5}
Assume that $(\A,\rho,\Sigma)$ is central, irreducible and satisfies condition (I).
Let ${\frak L}$ be the associated $\lambda$-graph system to $(\A,\rho, \Sigma)$.
If $(\A, \rho, \Sigma)$ is effective and ${\Cal O}_{\frak L}$ is simple, purely infinite, then
$\ARL$ is simple, purely infinite.
\endproclaim
\demo{Proof}
It suffices to show that 
for any nonzero $X \in \ARL$, ther exist 
$A,B \in \ARL$ such that 
$AXB=1$.
By the previous lemma there exists 
a projection $Q \in {\Cal D}_\rho$ 
and a nonzero positive element $Z \in {\Cal F}_{\rho,k}^l$
for some $k \le l$ such that
$\| Q X^* X Q - Z \| < \epsilon$.
We may assume that  $\| \E_\rho(X^* X ) \| =1$
so that 
$ 1 - \epsilon < \| Z \| < 1 + \epsilon.$
By Lemma 5.3, take an element 
$V \in \ARL$ such that 
$$
V Z V^* = 1, \qquad 
\| V \| < \frac{1}{\sqrt{\| Z \|}} + \epsilon 
< \frac{1}{\sqrt{1 -\epsilon}} + \epsilon.
$$
It follows that
$$
\| V Q X^* X Q V^* - 1 \| 
< \| V \| ^2 \| Q X^* X Q - Z \|
< (\frac{1}{\sqrt{1 -\epsilon}} + \epsilon)^2 \cdot \epsilon.
$$
We may take $\epsilon > 0$ small enough so that 
$
\| V Q X^* X Q V^* - 1 \|  < 1
$
and hence 
$V Q X^* X Q V^*$ is invertible in $\ARL$.
Thus we  complete the proof.
\qed
\enddemo
\heading 6. Tensor products of $C^*$-symbolic dynamical systems
\endheading

In this section, we will consider tensor products between 
 $C^*$-symbolic dynamical systems and finite families of automorphisms of unital $C^*$-algebras.
This construction yields 
interesting examples of $C^*$-symbolic dynamical systems 
beyond $\lambda$-graph systems, 
that will be studied in the following sections.
Throughout this section, 
we fix a unital $C^*$-algebra $\B$ 
and a finite family of automorphisms
$\alpha_i \in \Aut(\B), i=1,\dots,N$ of $\B$.
Tensor products $\otimes$ between $C^*$-algebras always mean 
the minimal $C^*$-tensor products $\otimes_{\min}$.
We set $\Sigma = \{ \alpha_1,\dots,\alpha_N\}$.
Consider a $C^*$-symbolic dynamical system
$(\A, \rho, \Sigma)$.

\proclaim{Proposition 6.1}
For $\alpha_i \in \Sigma, i=1,\dots,N$, 
define $\rho^{\Sigma\otimes}_{\alpha_i}\in \End(\B\otimes\A)$
by setting  
$$
\rho^{\Sigma\otimes}_{\alpha_i}(b\otimes a) = \alpha_i(b) \otimes \rho_{\alpha_i}(a)
\qquad \text{ for } b \in \B, a \in \A.  
$$
Then 
$(\B \otimes \A, \rho^{\Sigma\otimes}, \Sigma)$
becomes a $C^*$-symbolic dynamical system over $\Sigma$ such that the presented subshift $\Lambda_{\rho^{\Sigma\otimes}}$
is the same as the subshift $\Lambda_{\rho}$
presented by $(\A,\rho,\Sigma)$.
\endproclaim
\demo{Proof}
We will first prove that
$(\B \otimes \A, \rho^{\Sigma\otimes}, \Sigma)$
is a $C^*$-symbolic dynamical system.
Since 
$\{ \rho_{\alpha_i} \}_{i=1}^N$ is essential,
for $\epsilon >0$, 
there exist 
$
x_{i,j}, y_{i,j} \in \A,\, j=1,\dots,n(i), \,i=1,\dots,N
$ 
such that 
$$
\| \sum_{i=1}^{N} \sum_{j=1}^{n(i)} 
x_{i,j} \rho_{\alpha_i}(1) y_{i,j} - 1 \| < \epsilon
$$
so that we have
$$
\| \sum_{i=1}^N \sum_{j=1}^{n(i)}(1 \otimes x_{i,j}) 
(\rho^{\Sigma\otimes}_{\alpha_i}(1))( 1\otimes y_{i,j})
- 1 \| < \epsilon.
$$ 
Hence the closed ideal generated by 
$\{ \rho^{\Sigma\otimes}_{\alpha_i}(1): i=1,\dots,N \}$
is all of $\B \otimes \A$,
so that 
$\{ \rho^{\Sigma\otimes}_{\alpha_i} \}_{i=1}^N$ 
is essential.

Since $\{ \rho_{\alpha_i} \}_{i=1}^N$ is faithful on $\A$,
the homomorphism
$\xi_\rho:\A \longrightarrow \oplus_{i=1}^N \A_i,$
where $\A_i = \A, i=1,\dots,N$
defined by 
$
 \xi_\rho(a)=  \oplus_{i=1}^N  \rho_{\alpha_i} (a) 
$
 is injective.
Consider the homomorphisms:
$$
\align
\id_\B \otimes \xi_\rho & : 
b \otimes a \in \B \otimes \A 
\rightarrow b \otimes \xi_\rho(a) \in \B \otimes \xi_\rho(\A), \\
\oplus_{i=1}^N(\alpha_i \otimes \id) & :
(b_i \otimes a_i)_{i=1}^N \in \oplus_{i=1}^N (\B \otimes\A_i) 
\rightarrow (\alpha_i(b_i) \otimes a_i)_{i=1}^N \in \oplus_{i=1}^N (\B \otimes\A_i).
\endalign
$$
Since 
$\B \otimes \xi_\rho(\A)$ 
is a subalgebra of
$\B \otimes(\oplus_{i=1}^N \A_i) 
=\oplus_{i=1}^N (\B \otimes\A_i)$
and
both $\id_B \otimes \xi_\rho $ 
and $\oplus_{i=1}^N(\alpha_i \otimes \id)$ are isomorphisms,
the composition
$\oplus_{i=1}^N(\alpha_i \otimes \id) \circ (\id \otimes \xi_\rho)$
is isomorphic.
Hence 
$$
\oplus_{i=1}^N \rho^{\Sigma\otimes}_{\alpha_i} 
= \oplus_{i=1}^N (\alpha_i \otimes \rho_{\alpha_i}): 
\B \otimes\A \rightarrow \otimes_{i=1}^N (\B \otimes\A_i)
$$
is injective.
This implies that 
$\{ \rho^{\Sigma\otimes}_{\alpha_i} \}_{i=1}^N$ is faithful.

By the equality
$$
\rho^{\Sigma\otimes}_{\alpha_{i_n}}\circ \cdots \circ 
\rho^{\Sigma\otimes}_{\alpha_{i_1}} (1) 
= \rho_{\alpha_{i_n}}\circ \cdots \circ \rho_{\alpha_{i_1}}(1)
$$
for $\alpha_{i_1},\cdots,\alpha_{i_n}\in \Sigma$,
the presented subshifts 
$\Lambda_{\rho^{\Sigma\otimes}}$
and 
$\Lambda_{\rho}$
coincide.
\qed
\enddemo
We denote by 
$\Lambda$  the presented subshift
$\Lambda_{\rho}
(= \Lambda_{\rho^{\Sigma\otimes}})$.
Let $S_{\alpha_i}$ be the generating partial isometries of $\ARL$
satisfying 
$S_{\alpha_i}^* x S_{\alpha_i} = \rho_{\alpha_i}(x)$ 
for
$ x \in \A, i=1,\dots,N$,
and
 $\widetilde{S}_{\alpha_i}$ those of 
$\BMARL$
satisfying 
$\widetilde{S}_{\alpha_i}^* y \widetilde{S}_{\alpha_i} 
= \rho^{\Sigma\otimes}_{\alpha_i}(y)$ 
for
$ y \in \BMA, i=1,\dots,N$.
\proclaim{Proposition 6.2}
There exists a unital embedding $\tilde{\iota}$
of $\ARL$ into $\BMARL$ in a canonical way.
\endproclaim
\demo{Proof}
Define the injective homomorphism
$\iota:\A \rightarrow \BMA$ by setting
$\iota(a) = 1 \otimes a$ for $a \in A$.
Since the equality 
$\widetilde{S}_{\alpha_i}^* \iota(a) \widetilde{S}_{\alpha_i} 
= \iota(\rho_{\alpha_i}(a))$ 
for
$ a \in \A, i=1,\dots,N$
holds,
there exists a homomorphism
$\tilde{\iota}$ from $\ARL$ to $\BMARL$
satisfying
$\tilde{\iota}(a) = 1 \otimes a,
\tilde{\iota}(S_{\alpha_i}) = \widetilde{S}_{\alpha_i}$
for
$ a \in \A, i=1,\dots,N$
by the universality of $\ARL$. 
Let
$\E_\rho:\ARL \rightarrow {\Cal F}_{\rho}$
and
$\E_{\rho^{\Sigma\otimes}}:
\BMARL \rightarrow {\Cal F}_{\rho^{\Sigma\otimes}}
$
be the canonical conditional expectations respectively.
Define the $C^*$-subalgebras
${\Cal F}_{({\Bbb C}\otimes \A, \rho^{\Sigma\otimes})}
\subset \CMARL
$
of
$\BMARL$
by setting
$$
\align
\CMARL 
= & C^*(1 \otimes a, \widetilde{S}_{\alpha_i}: a \in \A, i=1,\dots,N),\\ 
{\Cal F}_{({\Bbb C}\otimes \A, \rho^{\Sigma\otimes})}
= & C^*(\widetilde{S}_{\mu}(1 \otimes a) \widetilde{S}_{\nu}^*:
 a \in \A,\, \mu, \nu \in \Lambda^*, |\mu| = |\nu |).
\endalign
$$ 
The diagrams
$$
\CD
\ARL @>\tilde{\iota}>> \CMARL @. \hookrightarrow \BMARL \\
@V{\E_\rho}VV  @VV{\E_{\rho^{\Sigma\otimes}}|_{\CMA}}V  @VV{\E_{\rho^{\Sigma\otimes}}}V \\ 
{\Cal F}_{\rho} @>\tilde{\iota}|_{{\Cal F}_{\rho}}>> 
{\Cal F}_{({\Bbb C}\otimes \A, \rho^{\Sigma\otimes})} @. \hookrightarrow 
{\Cal F}_{\rho^{\Sigma\otimes}}
\endCD
$$
are commutative.
Since $\iota:\A \rightarrow \CMA$ is isomorphic,
so is 
the restriction 
$
\tilde{\iota}|_{{\Cal F}_{\rho}}:
{\Cal F}_{\rho} 
\rightarrow
{\Cal F}_{({\Bbb C}\otimes \A, \rho^{\Sigma\otimes})} 
$
of 
${\Cal F}_{\rho}$.
One indeed sees that 
the condition
$S_\mu a S_\nu^* \ne 0$ 
for some 
$a \in \A, |\mu| = |\nu|$ 
implies 
$\widetilde{S}_\mu (1 \otimes a) \widetilde{S}_\nu^* \ne 0$
because of the equality
$
\iota(\rho_\mu(1) a \rho_\nu(1)) = 
\widetilde{S}_\mu^*\widetilde{S}_\mu (1 \otimes a) \widetilde{S}_\nu^*\widetilde{S}_\nu.
$
For 
$\sum_{\mu,\nu \in \Lambda^k}S_\mu a_{\mu,\nu}S_\nu^* 
\in {\Cal F}_{\rho}$,
suppose that 
$\tilde{\iota}(\sum_{\mu,\nu \in \Lambda^k}S_\mu a_{\mu,\nu}S_\nu^*)=0$.
It follows that for any $\xi, \eta \in \Lambda^k$,
$$
0 = \widetilde{S}_\xi^* 
(\sum_{\mu,\nu \in \Lambda^k}
\widetilde{S}_\mu (1 \otimes a_{\mu,\nu}) \widetilde{S}_\nu^*) \widetilde{S}_\eta 
= \widetilde{S}_\xi^* \widetilde{S}_\xi ( 1 \otimes a_{\xi,\eta})
\widetilde{S}_\eta^* \widetilde{S}_\eta
$$
so that 
$0=\widetilde{S}_\xi ( 1 \otimes a_{\xi,\eta})
\widetilde{S}_\eta^*,$
and hence 
$ S_\xi a_{\xi,\eta} S_\nu^* =0$.
This implies that 
$
\tilde{\iota}|_{{\Cal F}_{\rho^k}}:
{\Cal F}_{\rho}^k 
\rightarrow
{\Cal F}_{({\Bbb C}\otimes \A, \rho^{\Sigma\otimes})}^k 
$
is injective
and so is 
$
\tilde{\iota}|_{{\Cal F}_{\rho}}:
{\Cal F}_{\rho} 
\rightarrow
{\Cal F}_{({\Bbb C}\otimes \A, \rho^{\Sigma\otimes})}. 
$
Therefore by using a routine argument,
one concludes that 
$\tilde{\iota}:\ARL \rightarrow \CMARL$
is injective and hence isomorphic.
\qed
\enddemo
Let us prove that $(\BMA,\rho^{\Sigma\otimes},\Sigma)$
satisfies  condition (I)
if 
$(\A,\rho,\Sigma)$
satisfies condition (I). 
The result will be used in the following sections.
We set the $C^*$-subalgebras
$
{\Cal D}_{(\CMA, \rho^{\Sigma\otimes})} 
\subset
{\Cal D}_{\rho^{\Sigma\otimes}}
$
of
$
{\Cal F}_{\rho^{\Sigma\otimes}} 
$
by setting
$$
\align
{\Cal D}_{\rho^{\Sigma\otimes}} 
  = & C^*(\widetilde{S}_\mu x \widetilde{S}_{\mu}^* : 
          \mu \in \Lambda^*, x \in \BMA), \\
{\Cal D}_{(\CMA, \rho^{\Sigma\otimes})} 
  = & C^*(\widetilde{S}_\mu (1 \otimes a) \widetilde{S}_{\mu}^*: 
          \mu \in \Lambda^*, a \in \A). \\
\endalign
$$
We may identify 
the subalgebra $\DA$ of ${\Cal F}_{\rho}$ 
with the subalgebra 
${\Cal D}_{(\CMA, \rho^{\Sigma\otimes})}$
of 
${\Cal F}_{({\Bbb C}\otimes \A,\rho^{\Sigma\otimes})} $
through the map $\tilde{\iota}$ as in the preceding proposition.

Let $\varphi \in \B^*$ be a faithful state on $\B$.
It is well-known that there exists a faithful projection 
$\Theta_\varphi: \BMA \rightarrow \A$ of norm one satisfying
$\Theta_\varphi(b\otimes a)= \varphi(b)a$ for $b\otimes a \in \B\otimes\A$.
\proclaim{Lemma 6.3}
Let $\varphi \in \B^*$ be a faithful state on $\B$
satisfying $\varphi \circ \alpha_i = \varphi, i=1,\dots,N$.
The projection $\Theta_\varphi: \BMA \rightarrow \A$ 
of norm one can be extended to a projection of norm one 
$\Theta_{\Cal D}: \DBMA \rightarrow \DA$
such that 
$\Theta_{\Cal D}(x) = x$ for $ x\in \DA$.
\endproclaim
\demo{Proof}
For $k \in \Bbb N$, define the $C^*$-subalgebras 
${\Cal D}_{\rho}^k $ of ${\Cal D}_\rho$
and
$
{\Cal D}_{\rho^{\Sigma\otimes}}^k 
$
of
$
{\Cal D}_{\rho^{\Sigma\otimes}} 
$
by setting 
$$
\align
{\Cal D}_{\rho}^k 
  = & C^*(S_\mu a S_{\mu}^* : \mu \in \Lambda^k, a \in \A), \\
{\Cal D}_{\rho^{\Sigma\otimes}}^k 
  = & C^*(\widetilde{S}_\mu x \widetilde{S}_{\mu}^*: 
          \mu \in \Lambda^k, x \in \BMA). \\
\endalign
$$
For $x_\mu \in \BMA, \xi \in \Lambda^k$,
the identities
$$
\align
\Theta_\varphi( 
 \widetilde{S}_\xi^* (\sum_{\mu\in \Lambda^k}\widetilde{S}_\mu x_\mu 
\widetilde{S}_\mu^*)\widetilde{S}_\xi)
&= \Theta_\varphi( (1 \otimes \rho_\xi(1)) x_\xi (1 \otimes \rho_\xi(1)) \\
&=\rho_\xi(1) \Theta_\varphi(x_\xi ) \rho_\xi(1)
= S_\xi^*  S_\xi \Theta_\varphi(x_\xi ) S_\xi^*  S_\xi 
\endalign
$$
hold, 
so that the map defined by
$\Theta_{\Cal D}^k:\DBMAk \rightarrow \DAk$
$$
\Theta_{\Cal D}^k(\sum_{\mu\in \Lambda^k}
\widetilde{S}_\mu x_\mu \widetilde{S}_\mu^* )
= \sum_{\mu\in \Lambda^k} S_\mu \Theta_\varphi(x_\mu ) S_\mu^*.
$$
is well-defined for each $k \in \Zp$.
We will next see the restriction of 
$\Theta_{\Cal D}^{k+1}$ to $\DBMAk$ coincides with 
$\Theta_{\Cal D}^k$.
Since 
$
\sum_{\mu\in \Lambda^k}\widetilde{S}_\mu x_\mu \widetilde{S}_\mu^* 
\in \DBMAk
$
is written as 
$
\sum_{\mu\in \Lambda^k} \sum_{i=1}^N \widetilde{S}_\mu 
\widetilde{S}_{\alpha_i}\widetilde{S}_{\alpha_i}^* x_\mu 
\widetilde{S}_{\alpha_i}\widetilde{S}_{\alpha_i}^* \widetilde{S}_\mu^*
\in {\Cal D}_{\rho^{\Sigma\otimes}}^{k+1}, 
$
it follows that 
$$
\Theta_{\Cal D}^{k+1}(\widetilde{S}_\mu x_\mu \widetilde{S}_\mu^*)
 =
\sum_{i=1}^N \Theta_{\Cal D}^{k+1}( \widetilde{S}_{\mu \alpha_i} 
\rho^{\Sigma\otimes}_{\alpha_i}(x_\mu )
\widetilde{S}_{\mu \alpha_i}^*  
=
\sum_{i=1}^N  S_{\mu \alpha_i} 
\Theta_{\varphi}(\rho^{\Sigma\otimes}_{\alpha_i}(x_\mu ))
S_{\mu \alpha_i}^*. 
$$
As the state $\varphi$ is $\alpha_i$-invariant
for
$i=1,\dots,N$,
one has
for $\sum_j b_j \otimes a_j \in \BMA$,
$$
\align
\Theta_\varphi( \rho^{\Sigma\otimes}_{\alpha_i}(
\sum_j b_j \otimes a_j))
& = \sum_j \varphi(\alpha_i(b_j))\rho_{\alpha_i}(a_j)
=\sum_j \varphi(b_j)\rho_{\alpha_i}(a_j) \\
& =\rho_{\alpha_i}(\sum_j \varphi(b_j)a_j)
=S_{\alpha_i}^*\Theta_\varphi(\sum_j b_j\otimes a_j)S_{\alpha_i}
\endalign
$$
so that 
$
\Theta_\varphi( \rho^{\Sigma\otimes}_{\alpha_i}(x_\mu)) = 
S_{\alpha_i}^*\Theta_\varphi(x_\mu)S_{\alpha_i}
$
for $x_\mu \in \BMA$.
It then follows that 
$$
\Theta_{\Cal D}^{k+1}(\widetilde{S}_\mu x_\mu \widetilde{S}_\mu^*)
=
\sum_{i=1}^N  S_{\mu \alpha_i} 
S_{\alpha_i}^*\Theta_\varphi(x_\mu)S_{\alpha_i}
S_{\mu \alpha_i}^*
=S_\mu \Theta_\varphi(x_\mu) S_\mu^*
= \Theta_{\Cal D}^k(\widetilde{S}_\mu x_\mu \widetilde{S}_\mu^*). 
$$
Therefore the sequence
$\{ \Theta_{\Cal D}^k \}_{k=1}^{\infty}$
defines a projection from 
$\DBMA$ onto $\DA$, which we denote by $\Theta_{\Cal D}$.
\qed
\enddemo
\proclaim{Lemma 6.4}
Assume that $(\A,\rho,\Sigma)$ is central.
Then 
$\widetilde{S}_\mu (1 \otimes a ) \widetilde{S}_\mu^*$ 
commutes with $b \otimes 1$ for
$ a \in \A, \mu \in \Lambda^*$ and $b \in \B$.
\endproclaim
\demo{Proof}
Since
$
(1 \otimes a ) \rho^{\Sigma\otimes}_\mu (b \otimes 1)
= \rho^{\Sigma\otimes}_\mu (b \otimes 1) (1 \otimes a),
$
it follows that
$$
\align
  & \widetilde{S}_\mu (1 \otimes a ) \widetilde{S}_\mu^* (b \otimes 1)\\
= & \widetilde{S}_\mu (1 \otimes a ) \rho^{\Sigma\otimes}_\mu (b \otimes 1)\widetilde{S}_\mu^* 
 = \widetilde{S}_\mu \widetilde{S}_\mu^* (b \otimes 1) 
\widetilde{S}_\mu(1\otimes a) \widetilde{S}_\mu^*  
 = (b \otimes 1) 
\widetilde{S}_\mu(1\otimes a) \widetilde{S}_\mu^*.  \\
\endalign
$$
\qed
\enddemo
\proclaim{Theorem 6.5}
Assume that there exists a faithful state $\varphi$ on $\B$ 
invariant under $\alpha_i \in \Aut(\B), i=1,\dots,N$.
Suppose that $(\A,\rho,\Sigma)$ is central.
If $(\A,\rho,\Sigma)$ satisfies condition (I), 
then
$(\BMA, \rho^{\Sigma\otimes},\Sigma)$
satisfies condition (I) and is central.
\endproclaim
\demo{Proof}
Since $(\A,\rho,\Sigma)$ satisfies condition (I),
there exists a increasing sequence
$
\A_l, l\in \Zp
$
of $C^*$-subalgebras of $\A$ and  
a projection 
$q_k^l \in {\Cal D}_{\rho}\cap {\A_l}^\prime$
with $l \ge k$
satisfying the conditions of condition (I).
We set
$(\BMA)_l = \B \otimes\A_l, l\in \Zp.$
Then the conditions
$\overline{\cup_{l\in \Bbb N} (\BMA)_l} = \BMA$
and
$\rho^{\Sigma\otimes}_{\alpha_i}((\BMA)_l)
\subset (\BMA)_{l+1}$
are easy to veryfy.
Let $\tilde{\iota}:\ARL \hookrightarrow \BMARL$
be the embedding in Proposition 6.2. 
Put
$\tilde{q}_k^l = \tilde{\iota}(q_k^l) \in \DBMA$
for $l \ge k$.
By the preceding lemma, one sees that 
$\tilde{q}_k^l \in \DBMA \cap ((\BMA)_l)'$.
We will show that 
$\tilde{q}_k^l x \ne 0$ for $0 \ne x \in (\BMA)_l$.
As $xx^* \in \B\otimes\A_l$, one has
$\Theta_{\Cal D}(xx^*) = \Theta_{\varphi}(xx^*) \in \A_l$.
Hence
$q_k^l\Theta_{\varphi}(xx^*)\ne 0$.
By the equality
$
\Theta_{\Cal D}(\tilde{q}_k^l xx^*\tilde{q}_k^l) 
= q_k^l \Theta_{\varphi}(xx^*)q_k^l,
$
one obtains 
$\tilde{q}_k^l x\ne 0$.
Let 
$\tilde{\phi}_{\rho^{\Sigma\otimes}}(X) 
= \sum_{i=1}^N \widetilde{S}_{\alpha_i}X \widetilde{S}_{\alpha_i}^*
$
for
$X \in \DBMA$.
One has
$$
\tilde{q}_k^l \tilde{\phi}_{\rho^{\Sigma\otimes}}^m(\tilde{q}_k^l)
= \tilde{\iota}(q_k^l \phi_{\rho}^m(q_k^l) )= 0 
\qquad
\text{ for all } m= 1,2,\dots, k.
$$
Thus
$(\BMA, \rho^{\Sigma\otimes},\Sigma)$
satisfies condition (I).
If  $(\A,\rho,\Sigma)$ is central,
the projections $1 \otimes \rho_\mu(1)$ for $\mu \in \Lambda^*$
commute with $\BMA$,
so that 
$(\BMA, \rho^{\Sigma\otimes},\Sigma)$
is central.
\qed
\enddemo

We will study structure of the fixed point algebra
$\FRS$ 
of 
$\BMARL$
under the gauge action
$\widehat{\rho^{\Sigma\otimes}}$.
Recall that $\FR$ denote  the fixed point algebra of $\ARL$ 
under the gauge action $\hat{\rho}$.
Recall that for $k \in \Zp$ the $C^*$-subalgebras
$\FRK$ of $\FR$ and $\FRSK$ of $\FRS$ are  generated by
$S_\mu a S_\nu^*$ for $ \mu,\nu\in \Lambda^k, a \in \A$
and
$\widetilde{S}_\mu x \widetilde{S}_\nu^*$ 
for $ \mu,\nu\in \Lambda^k, x \in \BMA$
respectively.
Then we have
\proclaim{Lemma 6.6}
The map
$
\varPhi^k:\widetilde{S}_\mu(b \otimes a)\widetilde{S}_\nu^* \rightarrow 
b \otimes S_\mu a S_\nu^*
$
for $ b \otimes a \in \BMA, \ \mu, \nu\in \Lambda^k$
extends to an isomorphism from
$\FRSK$ to $\B \otimes \FRK$.
\endproclaim
\demo{Proof}
For 
$
Y = \sum_{\mu,\nu \in \Lambda^k}
\widetilde{S}_\mu( \sum_{j=1}^n b_j \otimes a_j)\widetilde{S}_\nu^*
\in \FRSK$,
put
$$
\varPhi^k(Y) = \sum_{j=1}^n 
(b_j \otimes \sum_{\mu,\nu \in \Lambda^k} S_\mu a_j S_\nu^*) 
\in \B \otimes \FRK.
$$
It follows that for $\xi,\eta \in \Lambda^k$
$$
\widetilde{S}_\xi^* Y \widetilde{S}_\eta
 =\widetilde{S}_\xi^* \widetilde{S}_\xi
( \sum_{j=1}^n b_j \otimes a_j)\widetilde{S}_\eta^*\widetilde{S}_\eta 
 = \sum_{j=1}^n b_j \otimes S_\xi^* S_\xi a_j S_\eta^* S_\eta 
 = (1 \otimes S_\xi^* ) \varPhi^k(Y) 
(1 \otimes S_\eta ) 
$$
Hence
$Y =0$ if and only if $\varPhi^k(Y)=0$.
As $\varPhi^k$ is a homomorphism from
$\FRSK$ to $\B \otimes \FRK$,
it yields an isomorphism.
\qed
\enddemo
The following lemma is straightforward.
\proclaim{Lemma 6.7}
Let $\alpha\otimes\iota^k_\rho: 
\B \otimes \FRK \rightarrow \B\otimes{\Cal F}_\rho^{k+1}$
be the homomorphism defined by 
$$
(\alpha\otimes\iota^k_\rho)( b \otimes S_\mu a S_\nu^*) 
= \sum_{i=1}^n \alpha_i(b) \otimes S_{\mu \alpha_i}\rho_{\alpha_i}(a) S_{\nu \alpha_i}^* \quad \text{ for } b\otimes a \in \BMA, \ \mu, \nu \in \Lambda^k.
$$
Then the diagram
$$
\CD
\FRSK @> \iota^k_{\rho^{\Sigma\otimes}} >> {\Cal F}_{\rho^{\Sigma\otimes}}^{k+1}\\
@V{ \varPhi^k }VV   @VV{ { \varPhi^{k+1} } } V\\
 \B \otimes \FRK @>>\alpha\otimes\iota^k_\rho > \B \otimes {\Cal F}_{\rho}^{k+1}\endCD
$$
is commutative, 
where 
$\iota^k_{\rho^{\Sigma\otimes}}:\FRSK \rightarrow {\Cal F}_{\rho^{\Sigma\otimes}}^{k+1}$ denotes the natural inclusion.
\endproclaim
Hence we have
\proclaim{Proposition 6.8}
The $C^*$-algebra
$
\FRS
$
is the inductive limit
$$
 \B \otimes {\Cal F}_{\rho}^{1}
\overset{\alpha\otimes\iota^1_{\rho}}\to{\longrightarrow}
   \B \otimes {\Cal F}_{\rho}^{2}
\overset{\alpha\otimes\iota^2_{\rho}}\to{\longrightarrow}
   \B \otimes {\Cal F}_{\rho}^{3}
\overset{\alpha\otimes\iota^3_{\rho}}\to{\longrightarrow}
   \cdots.
$$
\endproclaim
Let $\B = C(X)$ be
the commutative $C^*$-algebra of all continuous functins 
on a compact Hausdorff space $X$ with a finite family 
$h_1,\dots,h_N$ of homeomorphisms on $X$. 
Define
$\alpha_i \in \Aut(C(X)), i=1,\dots,N$
by
$\alpha_i(f)(t) = f(h_i(t))$ 
for 
$ f\in C(X), t \in X$.
Put $\Sigma = \{ \alpha_1,\dots,\alpha_N\}$
Take 
$(\A_{\frak L},\rho^{\frak L},\Sigma)$ 
for a $\lambda$-graph system $\frak L$ over $\Sigma$
as $(\A,\rho,\Sigma)$.
Then the above $C^*$-algebra 
$
\FRS 
$
is an AH-algebra.
If in particular 
$X={\Bbb T}$,
the algbera is an A${\Bbb T}$-algebra.  
We will study these examples in the following sections.

\heading 7.
 $C^*$-symbolic dynamical systems from homeomorphisms and graphs
\endheading
Let $h_1,\dots, h_N$ be a finite family of homemorphisms
on a compact Hausdorff space $X$.
Put $\Sigma = \{h_1,\dots, h_N\}$.
Let ${\Cal G}$ be a left-resolving finite labeled graph 
$(G,\lambda)$ over $\Sigma$
with underlying finite directed graph $G$
and labeling map
 $\lambda:E\rightarrow \Sigma$.
We denote by 
$G =(V,E)$,
where
$
V=\{v_1,\dots,v_{N_0}\}
$ 
is the finite set of its vertices
and
$
E= \{ e_1,\dots,e_{N_1}\}
$ is the finite set of its directed edges.
As in the begining of Section 2, 
we have a $C^*$-symbolic dynamical system 
$({\Cal A}_{\Cal G}, \rho^{\Cal G}, \Sigma)$.
Identify the homeomorphisms
$h_i$ with the induced automorphisms
$\alpha_i$ on $C(X)$.
By Proposition 6.1,
the tensor product
$(C(X)\otimes{\A_{\Cal G}}, {(\rho^{\Cal G})}^{\Sigma\otimes},\Sigma)$
of $C^*$-symbolic dynamical system is defined.
Put $X_i = X, i=1,\dots,N_0$
and
$$
\A_{{\Cal G},X} 
= C(X) \otimes \A_{\Cal G} 
 =C(\sqcup_{i=1}^{N_0}X_i),
 \qquad
\rho^{{\Cal G},X}
=    {(\rho^{\Cal G})}^{\Sigma\otimes}.  
$$
We will study the $C^*$-symbolic dynamical system
$(\A_{{\Cal G},X}, \rho^{{\Cal G},X},\Sigma)$.
Note that
the presented subshift 
$
\Lambda_{ \rho^{{\Cal G},X}}
$
is the sofic shift 
$
\Lambda_{\Cal G}
$ 
presented by the labeled graph $\Cal G$.

For $u,v \in V$,
let
$
H_n(u,v)
$
be the set 
$(f_1,\dots,f_n)$
of $n$-edges of the graph ${\Cal G}$ satisfying
$
s(f_1) =u, t(f_{i}) = s(f_{i+1}), i=1,\dots,n-1,
$
and
$t(f_n) = v.$
We set
$$
H_n(u) = \cup_{v \in V} H_n(u,v),\qquad
H_{\Cal G}^n = \cup_{u \in V} H_n(u), \qquad
H_{\Cal G} = \cup_{n=1}^\infty H_{\Cal G}^n.
$$
Then $\gamma =(f_1,\dots,f_n)\in H_n(v_i,v_j)$
yields a homeomorphism
$\lambda(\gamma)$ from $X_i$ to $X_j$ by setting
$$
\lambda(\gamma)(x) =
\lambda(f_n)\circ \cdots \circ\lambda(f_1)(x) 
\qquad
\text{ for } x \in X_i.
$$
For 
$x \in X_k$ with $k \ne i$, 
 $
\lambda(\gamma)(x) 
$ is not defined.
We set 
for $x \in X_i$
$$
orb_n(x) = 
\cup \{ \lambda(\gamma)(x) \mid \gamma \in H_n(v_i) \} 
\subset \sqcup_{j=1}^{N_0}X_j,
\qquad
orb(x) = \cup_{n=0}^{\infty}orb_n(x),
$$
where
$orb_0(x) =\{ x\}$.

\noindent
{\bf Definition.}
A family $(h_1,\dots,h_N)$ of homeomorphisms on $X$
is called ${\Cal G}$-{\it minimal}
if for any $ x\in \sqcup_{j=1}^{N_0}X_j$, the orbit
$orb(x)$ is dense in $\sqcup_{j=1}^{N_0}X_j$.  

\proclaim{Lemma 7.1}
The following conditions are equivalent:
\roster
\item"(i)" $(h_1,\dots,h_N)$ 
is  ${\Cal G}$-minimal;
\item"(ii)" There exists no proper closed subset 
$F \subset \sqcup_{j=1}^{N_0}X_j$
such that $\lambda(e_i)(F) \subset F$ for all $i=1,\dots,N_1$;
\item"(iii)" There exists no proper closed subset 
$F \subset \sqcup_{j=1}^{N_0}X_j$
such that $\cup_{i=1}^{N_1} \lambda(e_i)(F)= F$.
\endroster
\endproclaim
\demo{Proof}
(i)$ \Rightarrow$(ii)
If there exists a closed subset $F \subset \sqcup_{j=1}^{N_0}X_j$ such that 
$\lambda(e_i)(F) \subset F$ for all $i=1,\dots,N_1$,
take $x \in F \cap X_j$ for some $ j$.
Then $orb(x)$ is not dense in $ \sqcup_{j=1}^{N_0}X_j$. 

(ii)$ \Rightarrow$(i)
For  $x \in \sqcup_{j=1}^{N_0}X_j$,
let $F$ be the closure of $orb(x)$.
Then we have $\lambda(e_i)(F) \subset F$ for all $i=1,\dots,N_1$,
and
hence $F = \sqcup_{j=1}^{N_0}X_j$.

(ii)$ \Rightarrow$(iii) This implication is trivial.

(iii)$ \Rightarrow$(ii)
Suppose that there exists a closed subset $F \subset \sqcup_{j=1}^{N_0}X_j$ 
such that 
$\lambda(e_i)(F) \subset F$ for all $i=1,\dots,N_1$.
Put
$
\widetilde{F}_n 
= \cup_{\lambda(\gamma) \in H_{\Cal G}^n} \lambda(\gamma)(F)
$
a closed subset of $F$.
Since $\widetilde{F}_{n+1} \subset \widetilde{F}_n$ and   
$\sqcup_{j=1}^{N_0}X_j$ is compact, the set 
$E :=\cap_{n=1}^\infty \widetilde{F}_n$ is a nonempty closed
subset of  $\sqcup_{j=1}^{N_0}X_j$. 
Since 
$
\cup_{i=1}^{N_1}\lambda(e_i)(\widetilde{F}_n) = \widetilde{F}_{n+1},
$
one has 
$
\cup_{i=1}^N\lambda(e_i)(E) \subset E.
$
On the other hand, 
take $s(i)=1,\dots,N_0$ 
such that $v_{s(i)} = s(e_i)$.
Then we have
$$
\align  
& \cap_{n=1}^\infty \lambda(e_i)(\widetilde{F}_n)
= \cap_{n=1}^{\infty} \sqcup_{j=1}^{N_0} \lambda(e_i)( 
 \widetilde{F}_n \cap X_j)
=  \cap_{n=1}^{\infty} 
\lambda(e_i)( 
 \widetilde{F}_n \cap X_{s(i)})\\
\subset
&
\sqcup_{j=1}^{N_0} \cap_{n=1}^{\infty} 
\lambda(e_i)( 
 \widetilde{F}_n \cap X_j)
 =
\lambda(e_i)(E).
\endalign
$$
For
$
x \in \cap_{n=1}^\infty \cup_{i=1}^{N_1} \lambda(e_i)(\widetilde{F}_n) 
$
and
$n \in \Bbb N$,
 there exits $i_n = 1,\dots,N_1$ 
 such that 
$ x \in \lambda(e_{i_n})(\widetilde{F}_n).$
Find $i(x) =1,\dots,N_1$ such that 
$i(x)$ appears in 
$\{ i_n \mid n\in \Bbb N \}$ infinitely many times.
Since $\widetilde{F}_n, n\in \Bbb N$ are decreasing subsets,
one has 
$ x\in \lambda(e_{i(x)})(\widetilde{F}_n)$ for all $n \in \Bbb N$.
Hence 
$
x \in \cup_{i=1}^{N_1} \cap_{n=1}^\infty \lambda(e_i)(\widetilde{F}_n)
$
so that we have
$
\cup_{i=1}^{N_1} \cap_{n=1}^\infty \lambda(e_i)(\widetilde{F}_n) 
\supset 
 \cap_{n=1}^\infty \cup_{i=1}^{N_1} \lambda(e_i)(\widetilde{F}_n). 
$
Thus we have
$$
\cup_{i=1}^{N_1} \lambda(e_i)(E) 
\supset 
\cup_{i=1}^{N_1} \cap_{n=1}^\infty \lambda(e_i)(\widetilde{F}_n) 
\supset 
 \cap_{n=1}^\infty \cup_{i=1}^{N_1} \lambda(e_i)(\widetilde{F}_n) 
= 
 \cap_{n=1}^\infty \widetilde{F}_{n+1} 
 =E. 
$$
\qed
\enddemo
The following lemma is direct.
\proclaim{Lemma 7.2}
Let $J$ be an ideal of $\A_{{\Cal G},X}$.
Denote by $F \subset \sqcup_{j=1}^{N_0} X_j$ 
the closed subset such that 
$J = \{ f \in C(\sqcup_{j=1}^{N_0} X_j) \mid f(x) =0 \text{ for } x \in F \}$.
Then we have
\roster
\item"(i)" $J$ is a $\rho^{{\Cal G},X}$-invariant ideal of 
$\A_{{\Cal G},X}$ 
if and only if $\lambda(e_i)(F) \subset F$ for all $i=1,\dots,N_1$. 
\item"(ii)" $J$ is a saturated  ideal of $\A_{{\Cal G},X}$ 
if and only if $\cup_{i=1}^N \lambda(e_i)(F) \supset F$. 
\item"(iii)" $J$ is a $\rho^{{\Cal G},X}$-invariant saturated ideal of 
$\A_{{\Cal G},X}$ 
if and only if $\cup_{i=1}^N \lambda(e_i)(F) = F$. 
\endroster
\endproclaim
Hence we have
\proclaim{Lemma 7.3}
The following conditions are equivalent:
\roster
\item"(i)" $(h_1,\dots,h_N)$ 
is  ${\Cal G}$-minimal;
\item"(ii)" There exists no proper $\rho^{{\Cal G},X}$-invariant ideal of 
$\A_{{\Cal G},X}$;
\item"(iii)" There exists no proper $\rho^{{\Cal G},X}$-invariant saturated 
ideal of $\A_{{\Cal G},X}$.
\endroster
\endproclaim
A finite labeled graph ${\Cal G}$ is said to satisfy condition (I)
if for every vertex $v_i$ there exists  distinct paths with distinct labeled edges both of whose sourses and terminals are the vertex $v_i$. 
We denote by 
 ${\Cal O}_{{\Cal G}, h_1,\dots,h_N}$
the $C^*$-symbolic crossed product 
${\A_{{\Cal G},X}}\rtimes_{\rho^{{\Cal G},X}}\Lambda_{\Cal G}$
for the $C^*$-symbolic dynamical system
$(\A_{{\Cal G},X},\rho^{{\Cal G},X},\Sigma)$.
Assume that there exists a faithful $h_i$-invariant probability measure
on $X$.
\proclaim{Theorem 7.4}
Suppose that the labeled graph satisfies condition (I).
$(h_1,\dots,h_N)$ 
is  ${\Cal G}$-minimal
if and only if
the $C^*$-algebra
$
{\Cal O}_{{\Cal G},h_1,\dots,h_N} 
$
is simple.
\endproclaim
\demo{Proof}
Suppose that there exists a proper ideal 
$\I$ of 
${\Cal O}_{{\Cal G},h_1,\dots,h_N}$.
Since the labeled graph ${\Cal G}$ satisfies condition (I),
the $C^*$-symbolic dynamical system
$(\A_{\Cal G},\rho^{\Cal G},\Sigma)$ satisfies condition (I)
([Ma2;Section 4]),
so that 
$(\A_{{\Cal G},X},\rho^{{\Cal G},X},\Sigma)$
satisfies condition (I)
by Theorem 6.5.
Hence $J :=\I \cap \A_{{\Cal G},X}$ is a nonzero 
$\rho^{{\Cal G},X}$-invariant saturated
ideal of $\A_{{\Cal G},X}$.
If $J = \A_{{\Cal G},X}$, 
then $\A_{{\Cal G},X} \subset \I$
and
$S_\alpha^* S_\alpha \in \I$ so that $S_\alpha \in \I$.
Hence $\I =  {\Cal O}_{{\Cal G},h_1,\dots,h_N}$.
Therefore $J$ is not a proper ideal of $\A_{{\Cal G},X}$,
and by Lemma 7.3 
$(h_1,\dots,h_N)$ 
is not ${\Cal G}$-minimal.

Next suppose that 
$(h_1,\dots,h_N)$ 
is not ${\Cal G}$-minimal.
By Lemma 7.3, 
there exists a proper $\rho^{{\Cal G},X}$-invariant saturated ideal 
$J$ of $\A_{{\Cal G},X}$.
The ideal  $\I_J$ of ${\Cal O}_{{\Cal G},h_1,\dots,h_N}$
generated by $J$ 
satisfies $\I_J \cap \A_{{\Cal G},X} = J$ 
by Proposition 4.5.
Hence $\I_J$ is a proper ideal of 
${\Cal O}_{{\Cal G}, \gamma_1,\dots,\gamma_N}$.
\qed
\enddemo

In [KW;Corollary 33], 
Kajiwara-Watatani have proved a similar result for the $C^*$-algebras from circle bimodules.

For a vertex $u \in V$
put
$
H_n[u] = H_n(u,u).
$
Then we have
\proclaim{Proposition 7.5}
Suppose that ${\Cal G}$ satisfies condition (I) and is irreducible.
If there exists a path
$(f_1,\dots,f_n) \in H_n[v_i]$ for some 
vertex $v_i \in V$ and $n \in \Bbb N$ such that 
 the homeomorphism 
 $\lambda(f_n) \circ \cdots \circ \lambda(f_1)$
 on $X_i$ is minimal,
 then $(h_1,\dots,h_N)$ is ${\Cal G}$-minimal. 
\endproclaim
\demo{Proof}
Put
$\xi=(f_1, \dots,f_n)$.
Then 
$
\lambda(\xi)
$
is a minimal homeomorphism on $X_i$.
For  vertices $v_j,v_k \in V$, 
we may take paths 
$\gamma \in \cup_{m=1}^\infty H_m(v_i,v_j)$ 
and 
$\gamma' \in \cup_{m=1}^\infty H_m(v_k,v_i)$.
Since for any 
$x \in X_i$, the orbit 
$\cup_{l=0}^\infty \lambda(\xi)^l(x)$ is dense in $X_i$,
the set for any $y \in X_k$ 
$
\cup_{l=0}^\infty 
\lambda(\gamma) \circ \lambda(\xi)^l \circ \lambda(\gamma')(y)
$ 
is dense in $X_j$.
Thus
$(h_1,\dots,h_N)$ is ${\Cal G}$-minimal.
\qed
\enddemo
The above discussions may be generalized to a $\lambda$-graph system
with a family $\{ h_1,\dots,h_N\}$
of homeomorphisms of a compact Hausdorff space $X$. 
\heading 8. Irrational rotaton Cuntz-Kriger algebras
\endheading
Let $X $ be  the circle ${\Bbb T}$ in the complex plane.
Take  an arbitrary finite family of real numbers 
$\{ \theta_1,\dots,\theta_N\} $ with $\theta_i \in [0,1)$. 
Let 
$\alpha_i \in \Aut(C({\Bbb T}))$ 
be the automorphisms of $C({\Bbb T})$ 
defined by 
$
\alpha_i(f)(t) = f(e^{2 \pi \sqrt{-1}\theta_i}t), 
f\in C({\Bbb T}), t \in \Bbb T
$
for 
$i=1,\dots,N$.
Put
$\Sigma= \{ \alpha_1,\dots,\alpha_N\}$.
Let ${\Cal G}$ be a finite  directed labeled graph 
$(G,\lambda)$ over $\Sigma$
with underlying finite directed graph $G =(V,E)$
and left resolving labeling
$\lambda:E \rightarrow \Sigma$.  
We denote by $\{v_1,\dots,v_{N_0}\}$
the vertex set $V$.  
In [KW], Kajiwara-Watatani have studied the $C^*$-algebras constructed from 
circle correspondences.
Their situation is more general than ours.

Assume that each vertex of $V$ has both an incoming edge 
and an outgoing edge. 
Then we have a $C^*$-symbolic dynamical system
as in the preceding sections, which we denote by
$(\A_{{\Cal G},{\Bbb T}},\rho_{\theta_1,\dots,\theta_N},\Sigma)$.
Its $C^*$-symbolic crossed product is denoted by
$\OGN$.
Let
$A^{\Cal G}$ 
be the matrix for ${\Cal G}$ defined in (2.1).
\proclaim{Proposition 8.1}
The $C^*$-algebra $\OGN$ is the universal unital $C^*$-algebra
generated by 
$N$ partial isometries $S_i, i=1,\dots,N$
and $N_0$ partial unitaries $U_j, j=1,\dots,N_0$   
subject to the following relations:
$$
\align
& \sum_{m=1}^N S_m^* S_m =1,\qquad
 \sum_{j=1}^{N_0} U_j^* U_j =1,\qquad U_i^* U_i =U_i U_i^* \\
& U_i S_n 
= \sum_{j=1}^{N_0} A^{\Cal G}(i,\alpha_n,j) e^{2 \pi \sqrt{-1}\theta_n} S_n U_j,\\
& S_n S_n^* U_i =U_iS_n S_n^* \qquad \text{ for } i=1,\dots,N_0, \ n=1,\dots,N 
\endalign
$$
such that
$$
K_i(\OGN) 
= {\Bbb Z}^{N_0} / (1 - A_{\Cal G}){\Bbb Z}^{N_0} \oplus \Ker(1 - A_{\Cal G})
\quad i=0,1,
$$
where $A_{\Cal G}$ is the $N_0 \times N_0$ matrix defined by
$A_{\Cal G}(i,j) = \sum_{\alpha \in \Sigma}A^{\Cal G}(i,\alpha,j).$
\endproclaim
\demo{Proof}
It suffices to show the formulae of $K$-groups.
Since $K_i(\A_{{\Cal G},{\Bbb T}}) = {\Bbb Z}^{N_0}, i=0,1$,
by [Pim] (cf. [KPW])  the six term
exact sequence of K-theory:
$$
\CD
{\Bbb Z}^{N_0} @>\id - A_{\Cal G} >> 
{\Bbb Z}^{N_0} @>\id >> 
 K_0(\OGN) \\
@AAA @. @VVV \\
 K_1(\OGN)
@<<\id<
 {\Bbb Z}^{N_0} @<< \id - A_{\Cal G}< 
 {\Bbb Z}^{N_0}. \\ 
\endCD
$$
holds so that one has the short exact sequences for $i=0,1$
$$
0 
\longrightarrow 
{\Bbb Z}^{N_0} / (1 - A_{\Cal G}){\Bbb Z}^{N_0} 
\longrightarrow
K_i(\OGN)
\longrightarrow
 \Ker(1 - A_{\Cal G})
 \longrightarrow
 0.
 $$
 They split because $\Ker(1 - A_{\Cal G})$ is free 
 so that the desired formulae hold.
\qed
\enddemo
We denote by ${\Cal O}_{\Cal G}$ the $C^*$-algebra of the labeled graph 
${\Cal G}$.
It is isomorphic to a Cuntz-Krieger algebra (cf, [BP],[Ca],[Ma2],[Tom]). 
For $i,j=1,\dots,N_0$,
let $f_1,\dots,f_m$ be the set of edges in ${\Cal G}$
whose source is $v_i$ and terminal is $v_j$.
Then we set 
$
A^{{\Cal G}_\theta}(i,j) 
=e^{2 \pi \sqrt{-1}\theta_{k_1}} + \cdots + e^{2 \pi \sqrt{-1}\theta_{k_m}}
$  formal sums
for $\lambda(f_l) =\alpha_{k_l}, l=1,\dots,m$.
We have $N_0 \times N_0$ matrix $A^{{\Cal G}_\theta}$
with entries in formal sums of nonnegative real numbers.
\proclaim{Proposition 8.2}
Suppose that the labeled graph ${\Cal G}$ satisfies condition (I)
and is irreducible.
 If there exists $n \in \Bbb N$ and $i=1,\dots,N_0$ such that 
 the $(i,i)$-component $ (A^{{\Cal G}_\theta})^n(i,i)$
of the $n$-th power of the matrix
$A^{{\Cal G}_\theta}$
contains an irrational angle of rotation, 
then 
$(\alpha_1,\dots,\alpha_N)$ is ${\Cal G}$-minimal,
so that the $C^*$-algebra $\OGN$ is simple, purely infinite.
\endproclaim
\demo{Proof}
One knows that $(\alpha_1,\dots,\alpha_N)$ is ${\Cal G}$-minimal
by Proposition 7.5.
It is easy to see that 
$(\A_{\Cal G}, \rho_{{\theta_1},\dots,{\theta_N}},\Sigma)$
is effective.
As the algebra
${\Cal O}_{\Cal G}$ is purely infinite,
so is  $\OGN$ by Theorem 5.5.
\qed
\enddemo
We will study the structure of both the algebra $\OGN$
and the fixed point algebra 
${\Cal F}_{{\Cal G},\theta,\dots,\theta_N}$ of $\OGN$
under the gauge action. 
We denote by
${\Cal F}_{\Cal G}$
the fixed point algebra
of 
${\Cal O}_{\Cal G}$
under the gauge action.  
\proclaim{Proposition 8.3}
Assume that the labeled graph ${\Cal G}$
satisfies condition (I).
\roster
\item"(i)"
 $\OGN$ 
is isomorphic to the crossed product 
$
{\Cal O}_{\Cal G} \rtimes_{\gamma_{\theta_1,\dots,\theta_N}}{\Bbb Z}
$
of the Cuntz-Krieger algebra 
${\Cal O}_{\Cal G}$ of the labeled graph ${\Cal G}$ 
by an automorphisms 
$\gamma_{\theta_1,\dots,\theta_N}$
of 
${\Cal O}_{\Cal G}$.
\item"(ii)"
${\Cal F}_{{\Cal G},\theta,\dots,\theta_N}$
 is an A${\Bbb T}$-algebra, that is isomorphic to 
the crossed product 
$
{\Cal F}_{\Cal G} \rtimes_{\gamma_{\theta_1,\dots,\theta_N}}{\Bbb Z}
$
of the AF-algebra ${\Cal F}_{\Cal G}$
by the automorphism defined by the restriction of  
$\gamma_{\theta_1,\dots,\theta_N}$
to
${\Cal F}_{\Cal G}$. 
\endroster
\endproclaim
\demo{Proof}
(i)
Put $E_i = U_i^* U_i, i=1,\dots,{N_0}$.
The relations
$$
\sum_{j=1}^{N_0} E_j = 1,\qquad
 S_n^* E_i S_n = \sum_{j=1}^{N_0} A^{\Cal G}(i,\alpha_n, j) E_j
 $$
 hold
 for $n=1,\dots,N, i=1,\dots,{N_0}$.
 Hence the $C^*$-subalgebra 
 $C^*(S_n, E_i:n=1,\dots,N, i=1,\dots,{N_0})$ of $\OGN$ 
 generated by   
 $S_n, E_i:n=1,\dots,N, i=1,\dots,{N_0}$
 is isomorphic to the Cuntz-Krieger algebra 
 ${\Cal O}_{\Cal G}$ of the labeled graph ${\Cal G}$. 
Put
$U = \sum_{i=1}^{N_0} U_i$ a unitary.
It is straightforwad to see the following relations hold:
$$
US_nU^* = e^{2 \pi \sqrt{-1}\theta_n}S_n,\qquad
UE_i = E_i U = U_i,
$$
for $ n=1,\dots,N, i=1,\dots,{N_0}$.
Since the algebra $\OGN$ is generated by
$S_n, E_i$ for $ n=1,\dots,N, i=1,\dots,{N_0}$
and by putting 
$$
\gamma_{\theta_1,\dots,\theta_N}(S_n) = e^{2 \pi \sqrt{-1}\theta_n}S_n,
\qquad
\gamma_{\theta_1,\dots,\theta_N}(E_i) = E_i
$$
one sees that $\OGN$ is the crossed product 
of 
$C^*(S_n, E_i:n=1,\dots,N, i=1,\dots,{N_0})$
by the automorphism
$\gamma_{\theta_1,\dots,\theta_N}$.

(ii)
The AF-algebra ${\Cal F}_{\Cal G}$
is regarded as the $C^*$-subalgebra of
$\OGN$ generated by the elements of the form:
$S_\mu E_i S_\nu^*,  \mu, \nu \in \Lambda^*, |\mu| = |\nu |, i=1,\dots,N_0$.
Under the identification,
the algebra
${\Cal F}_{{\Cal G},\theta,\dots,\theta_N}$
is generated by 
${\Cal F}_{\Cal G}$ and the above unitary $U$.
By
$
\gamma_{\theta_1,\dots,\theta_N}(S_\mu E_i S_\nu^*) =
e^{2\pi \sqrt{-1}(\theta_{\mu_1}+\cdots + \theta_{\mu_k}
-\theta_{\nu_1}-\cdots-\theta_{\nu_k})}     
S_\mu E_i S_\nu^*
$
 for
 $
\mu = (\mu_1,\dots,\mu_k),\ \nu = (\nu_1,\dots,\nu_k) \in \Lambda^k,
$
one knows that
${\Cal F}_{{\Cal G},\theta,\dots,\theta_N}$
is isomorphic to the crossed product 
$
{\Cal F}_{\Cal G} \rtimes_{\gamma_{\theta_1,\dots,\theta_N}}{\Bbb Z}
$
of ${\Cal F}_{\Cal G}$
by   
$\gamma_{\theta_1,\dots,\theta_N}$.
By Proposition 6.8, one sees that
${\Cal F}_{{\Cal G},\theta,\dots,\theta_N}$
is an A${\Bbb T}$-algebra.
\qed
\enddemo
\heading 9. Irrational rotaton Cuntz algebras
\endheading
In this section,
we treat special cases of the previous section.
We consider a labeled graph of N-loops with single vertex.  
Let $A= C(\Bbb T)$ and 
$\Sigma = \{ 1, \dots, N\}, N>1$.
Take real numbers 
$\theta_1,\dots,\theta_N \in [0,1).$
Define 
$\alpha_i(f)(z) = f(e^{2\pi \sqrt{-1} \theta_i}z)
$
for
$f \in C({\Bbb T}), z \in \Bbb T.$
We have a $C^*$-symbolic dynamical system
$(C(\Bbb T),\alpha, \Sigma)$. 
Since $\alpha_i,i=1,\dots,N$ are automorphisms,
the associated subshift is the full shift 
$\Sigma^{\Bbb Z}.$
We denote by ${\Cal O}_{\theta_1,\dots,\theta_N}$ 
the $C^*$-symbolic crossed product 
$C(\Bbb T) \rtimes_{\alpha} \Sigma^{\Bbb Z}$.
As the algebra ${\Cal O}_{\theta_1,\dots,\theta_N}$ is the universal $C^*$-algebra generated by $N$ isometries $S_i, i=1,\dots,N$ and one unitary $U$ 
subject to the relations:
$$
\sum_{j=1}S_jS_j^* =1, \qquad
S_i^* S_i =1,\qquad
US_i = e^{2\pi \sqrt{-1} \theta_i}S_i U, \qquad i=1,\dots,N,
$$
it
is realized as the ordinary crossd product
$
{\Cal O}_N \rtimes_{\gamma_{\theta_1,\dots,\theta_N}}\Bbb Z
$ 
of the Cuntz algebra
${\Cal O}_N$ by the automorphism
$\gamma_{\theta_1,\dots,\theta_N}$  
 defined by 
 $
 \gamma_{\theta_1,\dots,\theta_N}(S_i) 
 = e^{2\pi \sqrt{-1} \theta_i} S_i$.
The K-groups are
$$
 K_0({\Cal O}_{\theta_1,\dots,\theta_N}) \cong
 K_1({\Cal O}_{\theta_1,\dots,\theta_N}) \cong {\Bbb Z}/(N-1){\Bbb Z}.
$$  
By Theorem 5.5 and Theorem 7.4, one sees
\proclaim{Proposition 9.1}
The $C^*$-algebra  ${\Cal O}_{\theta_1,\dots,\theta_N}$
is simple if and only if 
at least one of $\theta_1,\dots,\theta_N$ is irrational.
In this case,  ${\Cal O}_{\theta_1,\dots,\theta_N}$
is pure infinite.
\endproclaim
 
\noindent
{\bf Remark.}
The algebra 
${\Cal O}_{\theta_1,\dots,\theta_N}$
is the crossed product 
${\Cal O}_N \rtimes_{\gamma_{\theta_1,\dots,\theta_N}}{\Bbb Z}$ 
of the Cuntz algebra ${\Cal O}_N$
by the automorphism 
$
\gamma_{\theta_1,\dots,\theta_N}.
$
The condition that 
at least one of $\theta_1,\dots,\theta_N$ is irrational
is equivalent to the condition that
the automorphisms
$
(\gamma_{\theta_1,\dots,\theta_N})^n
$ are outer for all 
$n \in {\Bbb Z}, n\ne 0$.
Hence by
[Ki], 
the assertion for the simplicity of 
${\Cal O}_{\theta_1,\dots,\theta_N}$
in Proposition 9.1
holds.

We will study the fixed point algebra,
denoted by $\FN$,
of
${\Cal O}_{\theta_1,\dots,\theta_N}$
under the gauge action.
It is generated by elements 
of the form $S_\mu f S_\nu^*$ for 
$f \in C({\Bbb T}), |\mu | = |\nu |$.
Let $\FNk$ be the $C^*$-subalgebra of $\FN$ 
generated by elements of the form
$f \in C({\Bbb T}), |\mu | = |\nu |=k$.
The map 
$$
S_\mu f S_\nu^* \in \FNk 
\rightarrow 
f \otimes S_\mu S_\nu^* \in C({\Bbb T}) \otimes M_{N^k}
$$
yields an isomorphism between
$\FNk$ and $C({\Bbb T}) \otimes M_{N^k}$.
Then the natural inclusion
$\FNk \hookrightarrow {\Cal F}_{\theta_1,\dots,\theta_N}^{k+1}$
through the identity 
$
S_\mu f S_\nu^* = 
\sum_{i=1}^N S_{\mu i} \alpha_i(f) S_{\nu i}^*
$
corresponds to the inclusion
$$
\align
& f \otimes e_{i,j} \in C({\Bbb T}) \otimes M_{N^k}  \\
\hookrightarrow
& 
{\bmatrix
\alpha_1(f) \otimes e_{i,j} & & & 0\\
& \alpha_2(f) \otimes e_{i,j} & & \\
& & \ddots & \\
0 & & & \alpha_N(f) \otimes e_{i,j}
\endbmatrix
}
\in C({\Bbb T}) \otimes M_{N^{k+1}}.
\endalign
$$
For $\mu = (\mu_1,\dots,\mu_k) \in \Sigma^k$, 
we set
$\alpha_\mu = \alpha_{\mu_k}\circ \cdots \circ \alpha_{\mu_1}$.
Since $\FN$ is an inductive limit of the inclusions
$
\FNk \hookrightarrow {\Cal F}_{\theta_1,\dots,\theta_N}^{k+1},\
k=1,2,\dots
$
as in Proposition 6.8,
it is an A${\Bbb T}$-algebra.

\proclaim{Proposition 9.2}
The $C^*$-algebra
$\FN$ is simple
if and only if 
$\theta_i - \theta_j$ is irrational
for some  $i,j = 1,\dots,N$.
\endproclaim
\demo{Proof}
It is not difficult to prove the assertion directly by looking at the above inclusions
$\FNk \hookrightarrow {\Cal F}_{\theta_1,\dots,\theta_N}^{k+1}, k \in \Bbb N.$
The following argument is a shorter proof by using [Ki].
Let ${\Cal F}_N$ be the UHF-algebra of type $N^{\infty}$,
that is the fixed point algebra of ${\Cal O}_N$ by the gauge action.
By Proposition 8.3, 
$
{\Cal F}_{\theta_1,\dots,\theta_N}
$
is the crossed product
$
 {\Cal F}_N \rtimes_{\gamma_{\theta_1,\dots,\theta_N}}{\Bbb Z}
$ 
where
$\gamma_{\theta_1,\dots,\theta_N}(S_\mu S_\nu^*) =
e^{2\pi \sqrt{-1}
(\theta_{\mu_1}+\cdots +\theta_{\mu_k}
-\theta_{\nu_1}-\cdots-\theta_{\nu_k})}     
S_\mu S_\nu^*$
for
$\mu =(\mu_1,\dots,\mu_k), \nu=(\nu_1,\dots,\nu_k) \in \Sigma^k$.
Hence the automorphisms
$
\gamma_{\theta_1,\dots,\theta_N}
$
is the product type automorphism
$
\prod^{\otimes}\Ad(u_\theta)
=\Ad(u_\theta) \otimes \Ad(u_\theta)\otimes \cdots 
$
for the unitary
$
u_\theta =
\bmatrix
e^{2\pi\sqrt{-1}\theta_1}&        & 0 \\
                         & \ddots & \\
                       0 &        &  e^{2\pi\sqrt{-1}\theta_N}
\endbmatrix
$
in
$
 M_N({\Bbb C})
$
under the canonical identification between 
${\Cal F}_N$ and $M_N \otimes M_N \otimes \cdots $.
Then the condition that 
$\theta_i - \theta_j$ is irrational
for some $i,j = 1,\dots,N$ 
is equivalent to the condition that 
$(\Ad(u_\theta))^n \ne \id$ for all $n \in {\Bbb Z}, n \ne 0$.
In this case,
the product type automorphisms 
$
(\prod^{\otimes}\Ad(u_\theta))^n
$
are outer
for all $n \in {\Bbb Z},n \ne 0$.
Hence by [Ki], the assertion holds
\qed
\enddemo

For $\{\theta_1,\dots,\theta_N\}$ and $n \in \Bbb N$,
put 
$$
S_n(\theta_1,\dots,\theta_N) 
= \{ \theta_{i_1} + \cdots + \theta_{i_n} \mid i_1,\dots,i_n = 1,\dots,N\}.
$$
then the sequence 
$\{ S_n(\theta_1,\dots,\theta_N) \}_{n \in \Bbb N}$
of finite sets 
is said to be uniformly distributed in ${\Bbb T}$ ([Ki2])
if 
$$
\lim_{n \to \infty} \frac{1}{N^n} 
\sum_{i_1,\dots,i_n =1}^N 
f(e^{2\pi \sqrt{-1}(\theta_{i_1}+\cdots+\theta_{i_n})})
=
\int_{\Bbb T} f(t) dt
\qquad \text{ for all } f \in C({\Bbb T}).
$$
The 
following lemma is easy
\proclaim{Lemma 9.3}
$\{ S_n(\theta_1,\dots,\theta_N) \}_{n \in \Bbb N}$
is  uniformly distributed in ${\Bbb T}$ 
if and only if 
$\theta_i - \theta_j$ is irrational
for some  $i,j = 1,\dots,N$.
\endproclaim
\demo{Proof}
$\{ S_n(\theta_1,\dots,\theta_N) \}_{n \in \Bbb N}$
is  uniformly distributed in ${\Bbb T}$ 
if and only if 
$$
\lim_{n \to \infty} \frac{1}{N^n} 
\sum_{i_1,\dots,i_n =1}^N 
e^{2\pi \sqrt{-1} \ell (\theta_{i_1}+\cdots+\theta_{i_n})}
= 0
\qquad \text{ for all } \ell \in {\Bbb Z}, \ell \ne 0. \tag 9.1
$$
Since
$
\sum_{i_1,\dots,i_n =1}^N 
e^{2\pi \sqrt{-1} \ell (\theta_{i_1}+\cdots+\theta_{i_n})}
= 
(e^{2\pi \sqrt{-1} \ell \theta_1}
+ \cdots +
e^{2\pi \sqrt{-1} \ell \theta_N})^n,
$ 
the condition (9.1) holds 
if and only if 
$$
|e^{2\pi \sqrt{-1} \ell \theta_1}
+ \cdots +
e^{2\pi \sqrt{-1} \ell \theta_N}| < N
\qquad \text{ for all } \ell \in {\Bbb Z}, \ell \ne 0. \tag 9.2
$$ 
The condition (9.2) is equivalent to the condition that
 $\theta_i - \theta_j$ is irrational for some  $i,j = 1,\dots,N$.
\qed
\enddemo
Thereofre 
we have
\proclaim{Theorem 9.4}
For $\theta_1,\dots, \theta_N \in [0,1)$, the following conditions are equivalent:
\roster
\item"(i)"
$\theta_i - \theta_j$ is irrational for some $i,j = 1,\dots,N$.
\item"(ii)" 
$\FN$ is simple.
\item"(iii)" 
$\FN$ has real rank zero.
\endroster
\endproclaim
\demo{Proof}
The equivalence between (i) and (ii) follows from 
Proposition 9.2.
It suffices to show the equivalence between (i) and (iii).
Since
$$
\Sp(\undersetbrace{n}\to{u_{\theta}\otimes \cdots \otimes u_{\theta}})
=S_n(\theta_1,\dots,\theta_N)
$$
and
$\gamma_{\theta_1,\dots,\theta_N}$
is a product type automorphism on $\prod^{\otimes}Ad(u_\theta)$
on the UHF-algebra ${\Cal F}_N$,
by [Ki2;Lemma 5.2]
the crossed product 
${\Cal F}_N \rtimes_{\gamma_{\theta_1,\dots,\theta_N}}{\Bbb Z}$
has real rank zero if and only if
$S_n(\theta_1,\dots,\theta_N)
$
is uniformly distrubuted in ${\Bbb T}$.
\qed
\enddemo

We note that 
by [Ki;Lemma 5.2], the crossed product 
${\Cal F}_N \rtimes_{\gamma_{\theta_1,\dots,\theta_N}}{\Bbb Z}$
 has real rank zero
if and only if 
$\FN$ has a unique trace.

Consequently we obtain
\proclaim{Theorem 9.5}
For $\theta_1,\dots, \theta_N \in [0,1)$, 
suppose that there exist $i,j = 1,\dots,N$ such that
$\theta_i - \theta_j$ is irrational.
Then the $C^*$-algebra $\FN$ is a unital simple A${\Bbb T}$-algebra of real rank zero with a unique tracial state such that 
$$
K_0(\FN) \cong {\Bbb Z}[\frac{1}{N}],
\qquad
K_1(\FN) \cong {\Bbb Z}.
$$
Hence 
$\FN$ is the Bunce-Deddens algebra of type $N^{\infty}$.
\endproclaim
\demo{Proof}
Since $K_i(C({\Bbb T}\otimes M_{N^k})= {\Bbb Z}, i=0,1$
and
the homomorphisms in Proposition 6.8
yield
the $N$-multiplications on 
$K_0(C({\Bbb T}\otimes M_{N^k})= {\Bbb Z} \rightarrow 
K_0(C({\Bbb T}\otimes M_{N^{k+1}})= {\Bbb Z}
$
and
the identities  on 
$K_1(C({\Bbb T}\otimes M_{N^k})= {\Bbb Z} \rightarrow 
K_1(C({\Bbb T}\otimes M_{N^{k+1}})= {\Bbb Z},
$
we get the K-theory formulae by  Proposition 6.8.
The obtained isomorphism
from
$K_0(\FN)$ 
to 
$
{\Bbb Z}[\frac{1}{N}]
$
preserves their order and 
maps the unit $1$ of $\FN$ to $1$ in 
${\Bbb Z}[\frac{1}{N}]$.
Hence
$\FN$ is isomorphic to the Bunce-Deddens algebra of type $N^{\infty}$.
\qed
\enddemo


\bigskip

{\it e-mail} :kengo{\@}yokohama-cu.ac.jp

\bye